\documentclass[12pt, reqno]{amsart}
\usepackage[utf8]{inputenc}
\usepackage{amsmath, amsthm, amscd, amsfonts, amssymb, graphicx, xcolor}
\usepackage[bookmarksnumbered, colorlinks, plainpages]{hyperref}

\newtheorem{theorem}{Theorem}[section]

\theoremstyle{definition}
\newtheorem{definition}[theorem]{Definition}

\newtheorem{conjecture}[theorem]{Conjecture}

\theoremstyle{remark}

\numberwithin{equation}{section}
\usepackage{amsmath, amsthm, amssymb, amsfonts}
\usepackage{graphicx, xcolor}
\usepackage{hyperref}
\usepackage{cleveref}
\usepackage[bookmarksnumbered, colorlinks, plainpages]{hyperref}
\usepackage{amsmath, amsthm, amssymb, hyperref}
\usepackage{geometry}
\usepackage{amsmath}   
\usepackage{amsthm}    
\usepackage{hyperref}  
\usepackage{cleveref} 
\usepackage{longtable}
\usepackage{array}  
\usepackage{booktabs}  
\usepackage{lmodern}   
\usepackage{multicol}
\usepackage{geometry}
\usepackage{float}
\usepackage{algorithm}
\usepackage{algpseudocode}
\usepackage{ragged2e}

\begin{document}
\setcounter{page}{1}

\color{darkgray}{
\centerline{\today}

\centerline{}

\centerline{}


\title[RERPHNBQF]{Redefining Euler-Rabinowitsch Polynomials with Heegner Number Based Quadratic Formulation}

\author[Sudarshan Kumaresan, Shipra Kumari, Neha Mishra  ]{Sudarshan Kumaresan\textsuperscript{1}, Shipra Kumari\textsuperscript{2*} and Neha Mishra\textsuperscript{3*}}

\address{$^{1}$ Department of Electrical Engineering, Veermata Jijabai Technological Institute, Mumbai, India.}
\email{\textcolor[rgb]{0.00,0.00,0.84}{skumaresan\_b22@et.vjti.ac.in}}

\address{$^{2}$ Department of Mathematics, Veermata Jijabai Technological Institute, Mumbai, India.}
\email{\textcolor[rgb]{0.00,0.00,0.84}{skumari@hs.vjti.ac.in}}

\address{$^{3}$ Department of Electrical Engineering, Veermata Jijabai Technological Institute, Mumbai, India.}
\email{\textcolor[rgb]{0.00,0.00,0.84}{nmishra\_@el.vjti.ac.in}}

\begin{abstract}
This paper introduces a novel class of prime-generating quadratic polynomials defined by
\(
f_{Z,k,H}(n) = n^2 - (2Zk - 1)n + \frac{(2Zk - 1)^2 + H}{4}
\)
where \( Zk \in \mathbb{Z}_{\geq 0} \) and \( H \) belong to the set of Heegner numbers. This form is closely related to the Euler-Rabinowitsch polynomials given by through specific substitutions. The structure enables algebraic tuning for prime-rich outputs and provides a deeper insight into the impact of Heegner numbers on prime distribution. Using tools such as Bateman-Horn conjecture and prime-counting functions, we demonstrate that this family can be optimized to generate a high density of primes. This work offers new directions for research on analytic number theory and potential applications in cryptography and signal processing.
\par
\noindent \textit{Keywords.} class numbers; real quadratic fields; prime-generating polynomials; cryptography; signal processing
\newline
\noindent \textit{2020 Mathematics Subject Classification.} Primary 11B83, 11Y16; Secondary 11A41, 94A60
\end{abstract} 

\maketitle


\section{Introduction and preliminaries}

Euler-Rabinowitsch polynomials are a notable class of quadratic polynomials in number theory, known for their ability to generate primes under certain parameter constraints. The general form is given by the following \cite{Euler-Rabinowitschpolynomials}.
\[
F_{\Delta, q}(x) = qx^2 + (\alpha_\Delta - 1)qx + \frac{(\alpha_\Delta - 1)^2q - \Delta}{4q} \tag{1.1}\label{1.1}
\]
where \( \Delta \) is a discriminant associated with an imaginary quadratic field, \( q \in \mathbb{N} \) is a square-free positive integer divisor of \( \Delta \), and \( \alpha_\Delta \in \mathbb{Z} \) is an integer defined based on the properties of the divisibility. More precisely, let \( \Delta = 4D = \bar{\Delta}^2 \) be a discriminant with associated radicand \( D \), and let \( q \in \mathbb{N} \) be a square-free divisor of \( \Delta \). The parameter \( \alpha_\Delta \) is defined as:
\[
\alpha_\Delta =
\begin{cases}
1, & \text{if } 4q \mid \Delta \\
2, & \text{otherwise}
\end{cases}
\]
These parameters are intricately tied to class numbers, modular forms, and the arithmetic of quadratic fields, and have been extensively used in prime-generating contexts.

In this paper, we propose a new family of prime-rich polynomials defined by the following:
\[\boxed{
f_{Z,k,H}(n) = n^2 - (2Zk - 1)n + \frac{(2Zk - 1)^2 + H}{4}} \tag{1.2}\label{1.2}
\]
where \( Z, k \in \mathbb{Z}_{\geq 0} \) and \( H \) is a Heegner number. Heegner numbers are special integers \( H \in \{1, 2, 3, 7, 11, 19, 43, 67, 163\} \) for which the imaginary quadratic field \( \mathbb{Q}(\sqrt{-H}) \) has class number one \cite{mollin1996quadratic}. Our polynomial closely mirrors the Euler–Rabinowitsch form expressed in eq.(\ref{1.1}) with the substitutions \( q = 1 \), \( \Delta = -H \), \( \alpha_\Delta = 2Zk \) and \(x = -n\) are applied. This deliberate change, while algebraically simple, introduces a distinct behavior that separates our formulation from the classical Euler–Rabinowitsch model in both structure and performance.
This equation has the remarkable ability to generate a vast and diverse class of polynomials with a high density of prime output over specific integer intervals. By carefully choosing the parameters $Z$, $k$, and $H$, we can recover many of the well-known prime-generating polynomials from history as special cases of this universal form, including \(n^2 + n + 41 \) (Euler's prime polynomial), \(n^2 + n + 17 \) (Legendre's polynomial), and \( n^2 - 79n + 1601 \) (Ribenboim’s polynomial). 

This proves that our formulation not only generalizes previously known prime-rich polynomials but also unifies them under a single symbolic framework. In fact, with this equation alone, an infinite number of prime-rich polynomials can be systematically generated. This property opens new doors for understanding the structure and distribution of primes through a polynomial lens.

The cryptographic significance of \(f_{(Z, k, H)}(n)\) lies in its ability to produce large primes consistently with tunable control. In modern cryptosystems, especially in public-key protocols such as RSA, the generation of large primes is a cornerstone requirement. Our polynomial provides an efficient, structured and potentially deterministic approach to generating such primes with less computational overhead than random primality testing methods. 

To rigorously validate the prime-rich nature of this function, we performed a detailed analysis through theoretical observations, empirical tests, graphical comparisons, and algorithmic simulations. However, before delving into the analytical and proof-based sections, it is necessary to understand several important foundational concepts and theorems in number theory. These will serve as the basis for appreciating the full scope and depth of the polynomial’s behavior.

\begin{theorem}[Prime Number Theorem]\cite{prime}
Let $\pi(x)$ be the number of primes $\leq x$. Then $\pi(x) \sim \frac{x}{\log x}$ as $x \to \infty$, i.e., $\lim_{x \to \infty} \frac{\pi(x)\log x}{x} = 1$. 
\end{theorem}

\begin{theorem}[Fermat's Theorem on Non-existence of Prime-only Polynomials]
There exists no non-constant polynomial with integer coefficients that yields only primes for all integer inputs \cite{fermat}. 
\end{theorem}

\begin{conjecture}[Bateman-Horn Conjecture]
Given irreducible integer polynomials $f_1, \ldots, f_k$ satisfying the admissibility condition, the number of $n \leq x$ such that all $f_i(n)$ are prime is asymptotically:
\[
P(x) \sim \frac{C}{\prod \deg f_i} \int_2^x \frac{dt}{(\log t)^k} \tag{1.3}\label{1.3}
 \]
where $C$ is a product over primes adjusted for local obstructions \cite{batemanhorn}.
\end{conjecture}

\begin{definition}
Heegner numbers are integers $d$ such that the imaginary quadratic field $\mathbb{Q}(\sqrt{-d})$ has class number 1 \cite{heegnernumbers}. There are exactly 9 such $d$: $1, 2, 3, 7, 11, 19, 43, 67, 163$.
\end{definition}

\begin{definition}
An integer $a$ is a quadratic residue modulo $n$ if there exists $x$ such that $x^2 \equiv a \pmod{n}$. Otherwise, $a$ is non-residue \cite{quadraticresidues}.
\end{definition}

\begin{definition}
The Euler–Mascheroni constant $\gamma$ is defined as \cite{eulermascheroni}:
\[
\gamma = \lim_{n \to \infty} \left( \sum_{k=1}^{n} \frac{1}{k} - \log n \right) \approx 0.5772\ldots \tag{1.4}\label{1.4}
\]
\end{definition}

\begin{definition}
For $\Re(s) > 1$, the Riemann zeta function is given by \cite{reimannzetafunction}:
\[
\zeta(s) = \sum_{n=1}^\infty \frac{1}{n^s} = \prod_{p \text{ prime}} \left(1 - \frac{1}{p^s}\right)^{-1} \tag{1.5}\label{1.5} \]
\end{definition}

\begin{conjecture}[Riemann Hypothesis]
All nontrivial zeros of $\zeta(s)$ lie on the line $\text{Re}(s) = \frac{1}{2}$ in the complex plane \cite{riemannhypothesis}.
\end{conjecture}

\begin{definition}
The approximation:
\[
\sum_{p \leq x} \frac{1}{p} \sim \log \log x + \gamma \tag{1.6}\label{1.6}
\]
shows that the reciprocal sum over primes diverges slowly. This follows from the logarithm of the Euler product for $\zeta(s)$ and underlies many results in the prime number theory \cite{ApproximationfromEuler'sProductfortheRiemannZetaFunction}.
\end{definition}

\begin{definition}
A probabilistic test based on Fermat’s Little Theorem. For integer \( n > 2 \) and base \( 1 < a < n - 1 \), if \( a^{n-1} \not\equiv 1 \pmod{n} \), then \( n \) is composite. If the congruence holds, \( n \) is a probable prime to base \( a \), although some composites (Fermat pseudoprimes) may still pass \cite{fermatlittletheorem}.
\end{definition}

\begin{definition}
A strong probabilistic primality test. For odd \( n > 2 \), write \( n - 1 = 2^r d \) with \( d \) odd. If a base \( a \) fails both \( a^d \equiv 1 \pmod{n} \) and \( a^{2^j d} \equiv -1 \pmod{n} \) for all \( j < r \), then \( n \) is composite. Otherwise, \( n \) is a probable prime to base \( a \) \cite{millerrabin}.
\end{definition}

\begin{definition}
A probabilistic test using Euler's criterion and the Jacobi symbol. For odd \( n > 2 \), if \( a^{(n-1)/2} \not\equiv \left( \frac{a}{n} \right) \pmod{n} \), then \( n \) is composite. Otherwise, it is a probable prime to the base \( a \) \cite{solovayStrassen}.
\end{definition}

\begin{definition}
This test uses the Jacobi symbol to evaluate the probable primality. For odd \( n \) and \( a \) co-prime to \( n \), if \( a^{(n-1)/2} \equiv \left( \frac{a}{n} \right) \pmod{n} \), then \( n \) is a probable prime. Often used in conjunction with the Solovay-Strassen test \cite{jacobiprimality}.
\end{definition}

\begin{definition}
Based on Wilson’s Theorem. An integer \( n > 1 \) is prime if and only if \( (n - 1)! \equiv -1 \pmod{n} \). Although exact and deterministic, it is computationally impractical for large \( n \) \cite{lagrangeprimality}.
\end{definition}

\begin{definition}
States that \( n > 1 \) is prime if and only if \( (n - 1)! + 1 \equiv 0 \pmod{n} \). This test is theoretically sound but inefficient due to factorial growth, making it unsuitable for large numbers \cite{wilsonprimality}.
\end{definition}

\begin{definition}
RSA (Rivest–Shamir–Adleman) is a widely used public-key cryptosystem based on the hardness of factoring large composite numbers. It involves key generation, encryption, and decryption \cite{RSA}.

\noindent In key generation, two large primes \( p \) and \( q \) are chosen, and \( n = p \cdot q \). The totient \( \phi(n) = (p-1)(q-1) \) is calculated, and an exponent \( e \) is selected such that \( \gcd(e, \phi(n)) = 1 \). The private exponent \( d \) is the modular inverse of \( e \mod \phi(n) \), satisfying \( d \cdot e \equiv 1 \mod \phi(n) \). The public key is \( (e, n) \), and the private key is \( (d, n) \).

\noindent\[\text{Encryption of message \( m \) is done via:} \quad c \equiv m^e \mod n \tag{1.7}\label{1.7}
\] 
\noindent\[ \text{Decryption recovers \( m \) by:} \quad
m \equiv c^d \mod n \tag{1.8}\label{1.8}
\]

\noindent RSA is secure as long as factoring \( n \) remains computationally infeasible for large \( n \) (typically \( \geq 2048 \) bits).
\end{definition}

\begin{definition}
ECC is a public-key cryptosystem based on the algebraic structure of elliptic curves over finite fields. It offers security comparable to RSA with smaller keys \cite{ECC}.

\noindent An elliptic curve over \( \mathbb{F}_p \) is defined by:
\[
y^2 \equiv x^3 + ax + b \pmod{p} \tag{1.9}\label{1.9}
\]
with the non-singularity condition:
\[
4a^3 + 27b^2 \not\equiv 0 \pmod{p} \tag{1.10}\label{1.10}
\]

The set of points on the curve forms an abelian group. ECC security is based on the difficulty of the Elliptic Curve Discrete Logarithm Problem (ECDLP): given \( P \) and \( Q = kP \), determine \( k \). Applications include ECDH (key exchange), ECDSA (digital signatures), and efficient public-key encryption. A 256-bit ECC key offers security comparable to a 3072-bit RSA key.
\end{definition}

\begin{definition}
The Shor algorithm is a quantum algorithm that solves integer factorization in polynomial time, threatening RSA and ECC. It runs in \( \mathcal{O}((\log N)^3) \), using the quantum fourier transform for period finding \cite{Shor'salgorithm}.

Given \( f(x) = a^x \mod N \), it finds the period \( r \), and uses:
\[
\gcd(a^{r/2} \pm 1, N) \tag{1.11} \label{1.11}
\]
to factor \( N \). If \( r \) is even and \( a^{r/2} \not\equiv -1 \mod N \), then a non-trivial factor is obtained. Though large-scale quantum computers are not yet available, Shor's algorithm highlights the need for post-quantum cryptography.
\end{definition}

After developing a firm understanding of the fundamental concepts surrounding our prime-rich polynomial $f_{Z,k,H}(n)$, let us now move toward it's core function.

\section{Empirical Analysis of the Polynomial Family \texorpdfstring{\(f_{Z,k,H}(n)\)}{f(Z,k,H)(n)} for Varying and Fixed \texorpdfstring{\(k\)}{k}}

\noindent As expressed in eq.(\ref{1.2}) using the defined family of prime-rich polynomials:
\[
f_{Z,k,H}(n) = n^2 - (2Zk - 1)n + \frac{(2Zk - 1)^2 + H}{4}
\]
We first fix the values of $Z = 1$ and $H = 163$, where $163$ is known to be the largest Heegner number and has shown strong prime-yielding behavior compared to other Heegner numbers. Changing the parameter $k$ from $0$ to $100$, we generate a sequence of quadratic polynomials $f_k(n)$ and compute the number of primes each polynomial produces over the first \(100\) integer inputs $n = 0$ to $99$.

\noindent The following table.(\ref{table1}) displays the results of this empirical experiment:

\begin{longtable}{|c|c|c|c|}

\caption{Prime Count for \texorpdfstring{$f_k(n) = n^2 - (2k - 1)n + \frac{(2k - 1)^2 + 163}{4}$}{fk(n) = n^2 - (2k - 1)n + \frac{(2k - 1)^2 + 163}{4}} for \texorpdfstring{$0 \leq n \leq 99$}{0 \leq n \leq 99}}\\
\hline
\textbf{$k$} & \textbf{$f_k(n)$} & \textbf{Prime Count} & \textbf{Comment} \\
\hline
\endfirsthead

\hline
\textbf{$k$} & \textbf{$f_k(n)$} & \textbf{Prime Count} & \textbf{Comment} \\
\hline
\endhead
\hline
\endfoot

\hline
\endlastfoot
\label{table1}
\textbf{0} & \textbf{$n^2 - (-1)n + 41$} & \textbf{86} & \textbf{Euler's Polynomial} \\
1 & $n^2 - 1n + 41$ & 86 &	Symmetric to Euler’s Polynomial \\
2 & $n^2 - 3n + 43$ & 86 &Still matching high prime count \\
3 & $n^2 - 5n + 47$ & 86 & Slightly shifted, consistent performance
\\
4 & $n^2 - 7n + 53$ & 87 & Stable prime count range\\
5 & $n^2 - 9n + 61$ & 87 & Stable prime count range\\
6 & $n^2 - 11n + 71$ & 87 & Stable prime count range\\
7 & $n^2 - 13n + 83$ & 87 & Stable prime count range\\
8 & $n^2 - 15n + 97$ & 87 & Stable prime count range\\
9 & $n^2 - 17n + 113$ & 88 & Small increase begins\\
10 & $n^2 - 19n + 131$ & 88 & Small increase begins\\
11 & $n^2 - 21n + 151$ & 89 & Gradual rising trend\\
12 & $n^2 - 23n + 173$ & 89 & Gradual rising trend\\
13 & $n^2 - 25n + 197$ & 90 &Plateau starts\\
14 & $n^2 - 27n + 223$ & 90 &Plateau starts\\
15 & $n^2 - 29n + 251$ & 90 &Plateau starts\\
16 & $n^2 - 31n + 281$ & 91 &Small increase continues\\
17 & $n^2 - 33n + 313$ & 91 &Small increase continues\\
18 & $n^2 - 35n + 347$ & 92 &Prime-rich region approaching\\
19 & $n^2 - 37n + 383$ & 93 &High-prime region\\
20 & $n^2 - 39n + 421$ & 93 &High-prime region\\
21 & $n^2 - 41n + 461$ & 93 &High-prime region\\
22 & $n^2 - 43n + 503$ & 93 &High-prime region\\
23 & $n^2 - 45n + 547$ & 93 &High-prime region\\
24 & $n^2 - 47n + 593$ & 94 &Peak plateau begins\\
25 & $n^2 - 49n + 641$ & 94 &Peak plateau begins\\
26 & $n^2 - 51n + 691$ & 94 &Peak plateau begins\\
27 & $n^2 - 53n + 743$ & 94 &Peak plateau begins\\
28 & $n^2 - 55n + 797$ & 94 &Peak plateau begins\\
29 & $n^2 - 57n + 853$ & 94 &Peak plateau begins\\
30 & $n^2 - 59n + 911$ & 94 &Peak plateau begins\\
31 & $n^2 - 61n + 971$ & 94 &Peak plateau begins\\
32 & $n^2 - 63n + 1033$ & 94 &Plateau continues \\
33 & $n^2 - 65n + 1097$ & 94 &Plateau continues\\
34 & $n^2 - 67n + 1163$ & 94 &Plateau continues\\
\textbf{35} & \textbf{$n^2 - 69n + 1231$} & \textbf{95} & \textbf{Local prime maximum} \\
\textbf{36} & \textbf{$n^2 - 71n + 1301$} & \textbf{95} & \textbf{Local prime maximum}\\
\textbf{37} & \textbf{$n^2 - 73n + 1373$} & \textbf{95} & \textbf{Local prime maximum}\\
\textbf{38} & \textbf{$n^2 - 75n + 1447$} & \textbf{95} & \textbf{Local prime maximum}\\
\textbf{39} & \textbf{$n^2 - 77n + 1523$} & \textbf{95} & \textbf{Local prime maximum}\\
\textbf{40} & \textbf{$n^2 - 79n + 1601$} & \textbf{95} & \textbf{Ribenboim’s Polynomial} \\
41 & $n^2 - 81n + 1681$ & 94 & Slight decline begins\\
42 & $n^2 - 83n + 1763$ & 93 &Downward trend\\
43 & $n^2 - 85n + 1847$ & 93 &Downward trend\\
44 & $n^2 - 87n + 1933$ & 94 &Small rebound\\
45 & $n^2 - 89n + 2021$ & 93 &Mid-range plateau\\
46 & $n^2 - 91n + 2111$ & 93 &Mid-range plateau\\
47 & $n^2 - 93n + 2203$ & 93 &Mid-range plateau\\
48 & $n^2 - 95n + 2297$ & 93 &Mid-range plateau\\
49 & $n^2 - 97n + 2393$ & 93 &Continuation of flat trend\\
50 & $n^2 - 99n + 2491$ & 92 &Prime count dip begins\\
51 & $n^2 - 101n + 2591$ & 93 &Recovery zone\\
52 & $n^2 - 103n + 2693$ & 93 &Recovery zone\\
53 & $n^2 - 105n + 2797$ & 93 &Recovery zone\\
54 & $n^2 - 107n + 2903$ & 93 &Recovery zone\\
55 & $n^2 - 109n + 3011$ & 93 &Recovery zone\\
56 & $n^2 - 111n + 3121$ & 94 &Small peak\\
57 & $n^2 - 113n + 3233$ & 93 &	Leveling again\\
58 & $n^2 - 115n + 3347$ & 93 &	Leveling again\\
59 & $n^2 - 117n + 3463$ & 94 &Peak before another rise\\
\textbf{60} & \textbf{$n^2 - 119n + 3581$} & \textbf{95} & \textbf{Second local maximum}\\
\textbf{61} & \textbf{$n^2 - 121n + 3701$} & \textbf{95} & \textbf{Second local maximum} \\
\textbf{62} & \textbf{$n^2 - 123n + 3823$} & \textbf{95} & \textbf{Second local maximum}\\
\textbf{63} & \textbf{$n^2 - 125n + 3947$} & \textbf{95} & \textbf{Second local maximum}\\
\textbf{64} & \textbf{$n^2 - 127n + 4073$} & \textbf{95} & \textbf{Second local maximum}\\
\textbf{65} & \textbf{$n^2 - 129n + 4201$} & \textbf{95} & \textbf{Second local maximum}\\
66 & $n^2 - 131n + 4331$ & 94 &Broad high plateau\\
67 & $n^2 - 133n + 4463$ & 94 &Broad high plateau\\
68 & $n^2 - 135n + 4597$ & 94 &Broad high plateau\\
69 & $n^2 - 137n + 4733$ & 94 &Broad high plateau\\
70 & $n^2 - 139n + 4871$ & 94 &Broad high plateau\\
71 & $n^2 - 141n + 5011$ & 94 &Broad high plateau\\
72 & $n^2 - 143n + 5153$ & 94 &Broad high plateau\\
73 & $n^2 - 145n + 5297$ & 94 &Broad high plateau\\
74 & $n^2 - 147n + 5443$ & 94 &Broad high plateau\\
75 & $n^2 - 149n + 5591$ & 94 &Plateau continues\\
76 & $n^2 - 151n + 5741$ & 94 &Plateau continues\\
77 & $n^2 - 153n + 5893$ & 93 &Start of gradual decline\\
77 & $n^2 - 153n + 5893$ & 93 &Start of gradual decline\\
78 & $n^2 - 155n + 6047$ & 93 &Start of gradual decline\\
79 & $n^2 - 157n + 6203$ & 93 &Start of gradual decline\\
80 & $n^2 - 159n + 6361$ & 93 &Start of gradual decline\\
81 & $n^2 - 161n + 6521$ & 93 &Start of gradual decline\\
82 & $n^2 - 163n + 6683$ & 92 &	Small drop\\
83 & $n^2 - 165n + 6847$ & 91 &Continued decrease\\
84 & $n^2 - 167n + 7013$ & 91 &Continued decrease\\
85 & $n^2 - 169n + 7181$ & 90 &Drop towards 90s\\
86 & $n^2 - 171n + 7351$ & 90 &Drop towards 90s\\
87 & $n^2 - 173n + 7523$ & 90 &Drop towards 90s\\
88 & $n^2 - 175n + 7697$ & 89 &Decline steady\\
89 & $n^2 - 177n + 7873$ & 89 &Decline steady\\
90 & $n^2 - 179n + 8051$ & 88 &Gradual reduction\\
91 & $n^2 - 181n + 8231$ & 88 &Gradual reduction\\
92 & $n^2 - 183n + 8413$ & 87 &End of prime-rich region\\
93 & $n^2 - 185n + 8597$ & 87 &End of prime-rich region\\
94 & $n^2 - 187n + 8783$ & 87 &End of prime-rich region\\
95 & $n^2 - 189n + 8971$ & 87 &End of prime-rich region\\
96 & $n^2 - 191n + 9161$ & 87 &End of prime-rich region\\
97 & $n^2 - 193n + 9353$ & 86 &Final tail, same as k = 0\\
98 & $n^2 - 195n + 9547$ & 86 &Final tail, same as k = 0\\
99 & $n^2 - 197n + 9743$ & 86 &Final tail, same as k = 0\\
100 & $n^2 - 199n + 9941$ & 86 & Final tail, same as k = 0\\
\hline 
\end{longtable} 

\noindent From the above empirical data of prime counts generated by the polynomial family $f_k(n)$ in the range $0 \leq n < 100$ given in table.(\ref{table1}), it is observed that a total of \textbf{76 polynomials} yield a prime count from \textbf{90 to 95}, inclusive. Among these, \textbf{12 polynomials} attain the maximum observed count of \textbf{95 primes}, specifically for the values $k = 35$ to $40$ and $k = 60$ to $65$. These intervals appear to represent prime-rich regions within the family. Such clustering of high prime counts around specific values of $k$ suggests that the behavior is not random but structurally guided, likely influenced by the arithmetic properties of the quadratic polynomial or its discriminant. The presence of these bands reveals a deeper number-theoretic significance of the polynomial family.

A graphical analysis is presented to illustrate how the number of primes generated by the polynomial \( f_k(n) \) varies as \( k \) ranges from 0 to 100. Each point on the curve corresponds to the number of primes produced by the polynomial for a given value of \( k \), evaluated over the range \( 0 \leq n < 100 \).
\begin{figure}[H]
    \centering
    \includegraphics[width=1\linewidth]{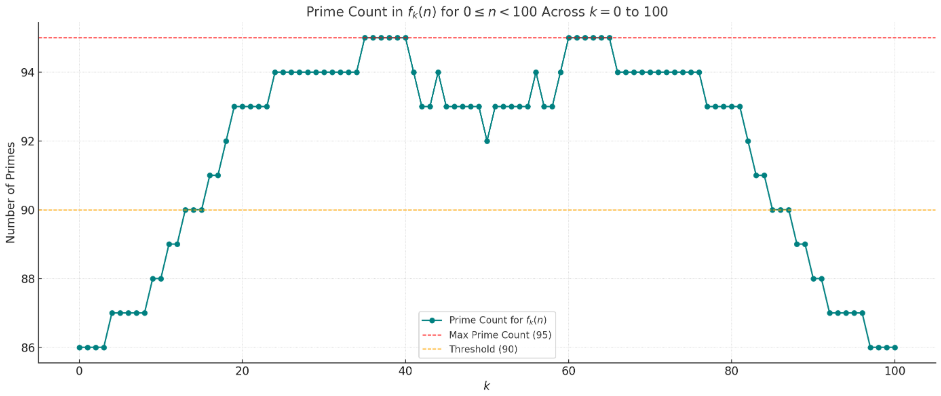}
    \caption{Graphical Analysis of Table.\ref{table1}}
    \label{g2}
\end{figure}

\noindent From fig.(\ref{g2}) the number of primes produced by the polynomial \( f_k(n) \) for each value of \( k \) in the range \( 0 \leq k \leq 100 \), we observe that out of \textbf{101} total polynomials, \textbf{68} lie above the threshold of 90 primes, \textbf{6} lie on the threshold while the remaining \textbf{27} lie below it. This indicates that a relatively large subset of the polynomial family demonstrates exceptionally strong prime-generating behavior. Notably, the distribution of these prime-rich polynomials is not random but exhibits a near-symmetric pattern around a central region (approximately \( k = 50 \)), suggesting an underlying arithmetic or algebraic symmetry in the structure of the family \( f_{Z,k,H}(n) \). This symmetry could be reflective of the balance between positive and negative coefficients of the polynomials and might point to deeper connections with known prime-rich forms such as Euler’s polynomial at \( k = 0 \) and Ribenboim’s polynomial at \( k = 40 \). This analysis supports the hypothesis that there exists a structured family of polynomials, governed by a parameter \( k \), that can consistently yield a high density of prime outputs.

To carry out a deeper analysis, we now fix the values of \( Z = 2 \), \( k = 40 \), and \( H = 163 \). This gives the polynomial \( f_{(2,40,163)}(n) = n^2 - 159n + 6361 \). We evaluate this polynomial for integer values of \( n \) from \(0\) to \(159\). The following table.(\ref{table2}) presents the value of \( n \), the result of the polynomial \( f_{(2,40,163)}(n) \), and whether that result is a prime number or not.\\

\vspace{1em}

\begin{longtable}{|c|c|c||c|c|c|}
\caption{Prime Counts for \texorpdfstring{$f_{(2, 40, 163)}(n) $}{f_{(2, 40, 163)}(n)} from \texorpdfstring{$ n = 0 $}{ n = 0 } to \texorpdfstring{$n = 159 $}{n = 159}}
\label{table2} \\
\hline
\textbf{n} & \textbf{\( f_{(2,40,163)}(n) \)} & \textbf{Status} & \textbf{n} & \textbf{\( f_{(2,40,163)}(n) \)} & \textbf{Status} \\
\hline
\endfirsthead

\hline
\textbf{n} & \textbf{\( f_{(2,40,163)}(n) \)} & \textbf{Status} & \textbf{n} & \textbf{\( f_{(2,40,163)}(n) \)} & \textbf{Status} \\
\hline
\endhead

\hline
\endfoot

\hline
\endlastfoot
 0 & 6361 & Prime & 80 & 41 & Prime \\
1 & 6203 & Prime & 81 & 43 & Prime \\
2 & 6047 & Prime & 82 & 47 & Prime \\
3 & 5893 & Not Prime & 83 & 53 & Prime \\
4 & 5741 & Prime & 84 & 61 & Prime \\
5 & 5591 & Prime & 85 & 71 & Prime \\
6 & 5443 & Prime & 86 & 83 & Prime \\
7 & 5297 & Prime & 87 & 97 & Prime \\
8 & 5153 & Prime & 88 & 113 & Prime \\
9 & 5011 & Prime & 89 & 131 & Prime \\
10 & 4871 & Prime & 90 & 151 & Prime \\
11 & 4733 & Prime & 91 & 173 & Prime \\
12 & 4597 & Prime & 92 & 197 & Prime \\
13 & 4463 & Prime & 93 & 223 & Prime \\
14 & 4331 & Not Prime & 94 & 251 & Prime \\
15 & 4201 & Prime & 95 & 281 & Prime \\
16 & 4073 & Prime & 96 & 313 & Prime \\
17 & 3947 & Prime & 97 & 347 & Prime \\
18 & 3823 & Prime & 98 & 383 & Prime \\
19 & 3701 & Prime & 99 & 421 & Prime \\
20 & 3581 & Prime & 100 & 461 & Prime \\
21 & 3463 & Prime & 101 & 503 & Prime \\
22 & 3347 & Prime & 102 & 547 & Prime \\
23 & 3233 & Not Prime & 103 & 593 & Prime \\
24 & 3121 & Prime & 104 & 641 & Prime \\
25 & 3011 & Prime & 105 & 691 & Prime \\
26 & 2903 & Prime & 106 & 743 & Prime \\
27 & 2797 & Prime & 107 & 797 & Prime \\
28 & 2693 & Prime & 108 & 853 & Prime \\
29 & 2591 & Prime & 109 & 911 & Prime \\
30 & 2491 & Not Prime & 110 & 971 & Prime \\
31 & 2393 & Prime & 111 & 1033 & Prime \\
32 & 2297 & Prime & 112 & 1097 & Prime \\
33 & 2203 & Prime & 113 & 1163 & Prime \\
34 & 2111 & Prime & 114 & 1231 & Prime \\
35 & 2021 & Not Prime & 115 & 1301 & Prime \\
36 & 1933 & Prime & 116 & 1373 & Prime \\
37 & 1847 & Prime & 117 & 1447 & Prime \\
38 & 1763 & Not Prime & 118 & 1523 & Prime \\
39 & 1681 & Not Prime & 119 & 1601 & Prime \\
40 & 1601 & Prime & 120 & 1681 & Not Prime \\
41 & 1523 & Prime & 121 & 1763 & Not Prime \\
42 & 1447 & Prime & 122 & 1847 & Prime \\
43 & 1373 & Prime & 123 & 1933 & Prime \\
44 & 1301 & Prime & 124 & 2021 & Not Prime \\
45 & 1231 & Prime & 125 & 2111 & Prime \\
46 & 1163 & Prime & 126 & 2203 & Prime \\
47 & 1097 & Prime & 127 & 2297 & Prime \\
48 & 1033 & Prime & 128 & 2393 & Prime \\
49 & 971 & Prime & 129 & 2491 & Not Prime \\
50 & 911 & Prime & 130 & 2591 & Prime \\
51 & 853 & Prime & 131 & 2693 & Prime \\
52 & 797 & Prime & 132 & 2797 & Prime \\
53 & 743 & Prime & 133 & 2903 & Prime \\
54 & 691 & Prime & 134 & 3011 & Prime \\
55 & 641 & Prime & 135 & 3121 & Prime \\
56 & 593 & Prime & 136 & 3233 & Not Prime \\
57 & 547 & Prime & 137 & 3347 & Prime \\
58 & 503 & Prime & 138 & 3463 & Prime \\
59 & 461 & Prime & 139 & 3581 & Prime \\
60 & 421 & Prime & 140 & 3701 & Prime \\
61 & 383 & Prime & 141 & 3823 & Prime \\
62 & 347 & Prime & 142 & 3947 & Prime \\
63 & 313 & Prime & 143 & 4073 & Prime \\
64 & 281 & Prime & 144 & 4201 & Prime \\
65 & 251 & Prime & 145 & 4331 & Not Prime \\
66 & 223 & Prime & 146 & 4463 & Prime \\
67 & 197 & Prime & 147 & 4597 & Prime \\
68 & 173 & Prime & 148 & 4733 & Prime \\
69 & 151 & Prime & 149 & 4871 & Prime \\
70 & 131 & Prime & 150 & 5011 & Prime \\
71 & 113 & Prime & 151 & 5153 & Prime \\
72 & 97 & Prime & 152 & 5297 & Prime \\
73 & 83 & Prime & 153 & 5443 & Prime \\
74 & 71 & Prime & 154 & 5591 & Prime \\
75 & 61 & Prime & 155 & 5741 & Prime \\
76 & 53 & Prime & 156 & 5893 & Not Prime \\
77 & 47 & Prime & 157 & 6047 & Prime \\
78 & 43 & Prime & 158 & 6203 & Prime \\
79 & 41 & Prime & 159 & 6361 & Prime \\
\hline
\hline
\end{longtable}
\begin{figure}[H]
    \centering
    \includegraphics[width=1\linewidth]{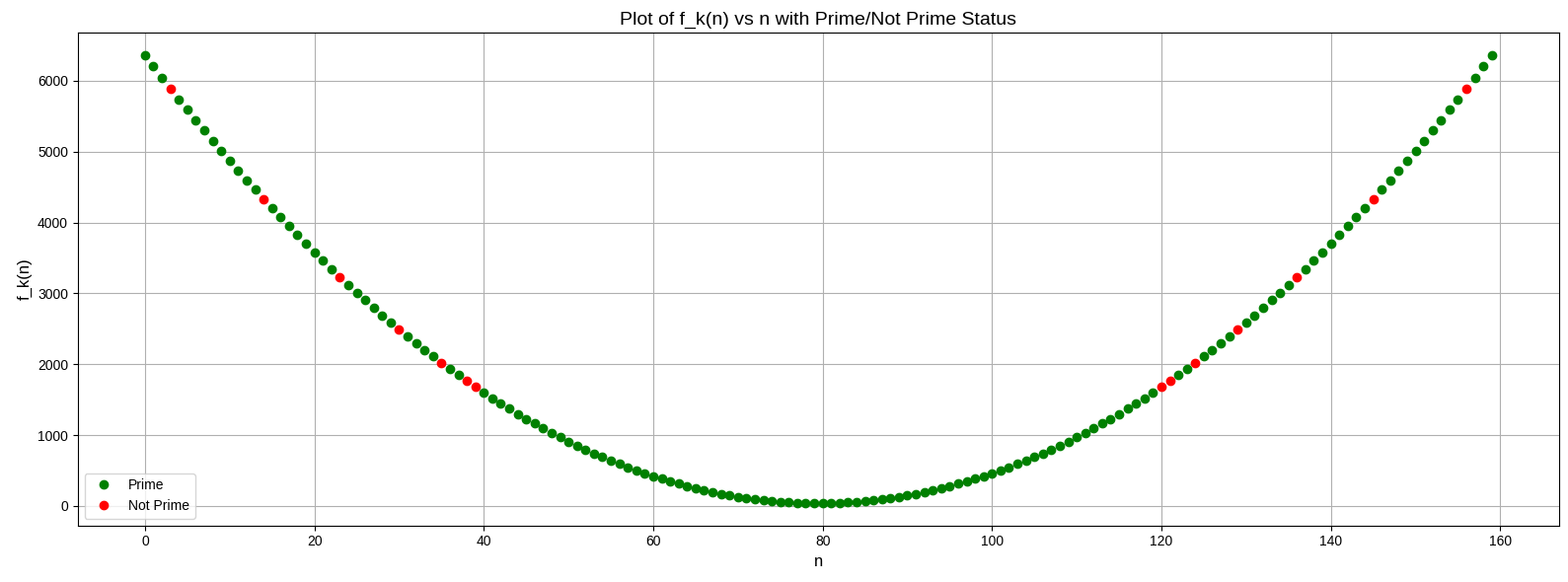}
    \caption{Graphical Analysis of Table.\ref{table2}}
    \label{f3}
\end{figure}
\paragraph{}
Table.(\ref{table2}) and fig.(\ref{f3}) showcase the output of the polynomial \( f_{(2,40,163)}(n) \) for \( n = 0 \) to \( n = 159 \), where we observe that out of \textbf{160} evaluated values, \textbf{146} are prime numbers and only \textbf{14} are not. This reflects a remarkably high prime density, suggesting that the polynomial under this configuration is strongly biased towards generating prime numbers. An intriguing pattern of symmetry is also evident in the data: the values of \( f(n) \) from \( n = 80 \) to \( n = 159 \) form a mirror image of the values from \( n = 0 \) to \( n = 79 \). Completing the square yields \( f(n) = (n - \frac{(2Zk-1)}{2})^2  + \frac{H}{4} \), indicating that the function is symmetric about \( n = \frac{(2Zk-1)}{2}\). When this expression is evaluated over a symmetric interval centered around \( n = 79.5 \), it results in a symmetric distribution of values. This mathematical symmetry is visually and numerically confirmed in the table. \textbf{Moreover, such symmetry is consistently observed for all natural values of \( k \) (except for \(k=0\)),} making it a general and intrinsic property of the polynomial family \( f_{(Z,k,H)}(n) \).

\paragraph{}
From the empirical analysis of the polynomial family \( f_{(Z,k,H)}(n) \), two key observations were made. First, when the parameter \( k \) is varied over a range of values while keeping \( Z \) and \( H \) constant, a clear symmetry emerges in the distribution of prime-rich polynomials. Specifically, the number of primes generated by the polynomial exhibits a mirrored pattern across certain values of \( k \), indicating that the behavior of prime generation is not random but follows an underlying structure dependent on the value of \( k \). 

\paragraph{}
Second, when \( k \) is fixed and \( n \) is varied over a given range, we observe a different form of symmetry; the actual values obtained from evaluating \( f_{(Z,k,H)}(n) \) display a mirror-like pattern about the midpoint of the \( n \) range. For example, when \( n \) varies from 0 to 159 with \( k = 40 \), the values of \( f_{(Z,k,H)}(n) \) from \( n = 0 \) to \( 79 \) are the mirror image of those from \( n = 80 \) to \( 159 \), resulting in a symmetric distribution of the outputs.

\paragraph{}
These two forms of symmetry highlight the distinct structural properties of the polynomial family. The symmetry observed when varying \( k \) reflects how the configuration of the polynomial parameters affects its overall prime-generating potential. In contrast, the symmetry observed when fixing \( k \) reveals that the polynomial itself inherently possesses a reflective nature with respect to \( n \), which is a consequence of its quadratic form. Together, these observations suggest that the family \( f_{(Z,k,H)}(n) \) is not only capable of generating primes densely in certain configurations but also encodes mathematical symmetry that can be exploited for further theoretical and cryptographic applications.

\section{Prime Rich Property Determination of \texorpdfstring{\(f_{(Z,k,H)}(n)\)}{f(Z,k,H)(n)} Using the Bateman-Horn Conjecture}

\noindent Let us denote eq.(\ref{1.2}) by:
\[
f(n) = n^2 - An + B,
\]
where \( A = (2Zk - 1) \) and \( B = \frac{A^2 + H}{4} \). 

\noindent So, the polynomial becomes:
\[
f(n) = n^2 - An + \frac{A^2 + H}{4}
\]

\noindent For a single irreducible polynomial \( f(n) \in \mathbb{Z}[n] \), the Bateman-Horn conjecture states that the number of integers \( n \in [1, N) \) such that \( f(n) \) is prime is asymptotically given by:
\[
\pi_f(N) \sim C_f \int_2^N \frac{dt}{\log t} \tag{3.1}\label{3.1}
\]
\noindent where \( C_f \) is the Bateman-Horn constant for the polynomial \( f(n) \), defined as:
\[
C_f = \prod_{\substack{p\ \text{prime}}} \left( 1 - \frac{\omega_f(p)}{p} \right) \left(1 - \frac{1}{p}\right)^{-1} \tag{3.2} \label{3.2}
\]
and \( \omega_f(p) \) is the number of distinct roots of \( f(n) \mod p \), i.e., the number of solutions to \( f(n) \equiv 0 \pmod{p} \)\cite{batemanhorn}.

\noindent Thus, the expected number of primes generated by \( f_{(Z,k,H)}(n) \) as \( n \in [0, N) \) is:
\[
\pi_{(Z,k,H)}(N) \sim C_{(Z,k,H)} \int_2^N \frac{dt}{\log t} \tag{3.3}\label{3.3}
\]

\noindent Applying prime number theorem, the integral is asymptotically equivalent to:
\[
\int_2^N \frac{dt}{\log t} = \operatorname{Li}(N) \sim \frac{N}{\log N} \tag{3.4}\label{3.4} \]

\noindent Hence, from eq.(\ref{3.3}) and eq.(\ref{3.4}) we obtain:
\[ \boxed{
\pi_{(Z,k,H)}(N) \sim C_{(Z,k,H)} \cdot \frac{N}{\log N}} \tag{3.5}\label{3.5}
\]

\noindent The \textbf{prime-rich nature} of the polynomial family is inferred when the constant \( C_{(Z,k,H)} \) is relatively large compared to other polynomials of similar degree.\\

\noindent To support this:
\begin{itemize}
    \item The discriminant of the polynomial \( f(n) \) is:
    \[
    \Delta = A^2 - 4B = A^2 - (A^2 + H) = -H,
    \]
    indicating that the polynomial is irreducible for positive \( H \), which satisfies the requirement of the conjecture.
    \item The number of roots modulo primes \( \rho_f(p) \) is small for large primes, contributing to a higher value of \( C_f \).
    \item When specific values of \( Z, k, H \) such as \( Z=2, k=40, H=163 \) are chosen, empirical evidence shows a high density of primes in the output values of \( f(n) \).
\end{itemize}

\noindent Eq.(\ref{3.5}) justifies the prime-rich behavior of the polynomial family \( f_{(Z,k,H)}(n) \), particularly for values like \( f_{(2,40,163)}(n) \), which produce an unusually high number of primes for small \( N \), consistent with a relatively large constant \( C_{(Z,k,H)} \). After analyzing the behavior of the polynomials belonging to the family \( f_{(Z,k,H)}(n) \) it was observed that these polynomials generated more primes in fact than expected using Bate-Horn conjecture.

\noindent In order to have a clearer understanding about the \textbf{prime-rich property} of the polynomial \( f_{(Z,k,H)}(n) \), let us take some examples where we compare the expected number of primes using the Bateman-Horn conjecture and the actual number of primes generated by the polynomial for a given range of \( n \).

\vspace{0.5em}
\noindent\textbf{Example 1.} Consider the values \( Z = 1 \), \( H = 163 \), and \( k = 100 \), with \( n \in [0, 199] \). The corresponding polynomial becomes:
\[
f_{(Z,k,H)}(n) = n^2 - (2Zk - 1)n + \frac{(2Zk - 1)^2 + H}{4} = n^2 - 199n + 9941.
\]

\noindent This is a quadratic polynomial of the form \( f(n) = n^2 + an + b \), where:
\[
a = -199, \quad b = 9941.
\]

\noindent The discriminant of the polynomial is:
\[
D = a^2 - 4b = (-199)^2 - 4 \cdot 9941 = 39601 - 39764 = -163.
\]

\noindent To compute the Bateman-Horn constant \( C_f \), we use the general formula:
\[
C_f = \prod_{p} \left(1 - \frac{\omega_f(p)}{p} \right) \left(1 - \frac{1}{p} \right)^{-1},
\]
\noindent where \( \omega_f(p) \) is the number of distinct solutions to \( f(n) \equiv 0 \pmod{p} \). Since \( f(n) \) is a quadratic, this depends on whether the discriminant \( D \) is a quadratic residue modulo \( p \), and hence:
\[
\omega_f(p) = 1 + \left( \frac{D}{p} \right) = 1 + \left( \frac{-163}{p} \right),
\]
\noindent where \( \left( \frac{\cdot}{p} \right) \) denotes the Legendre symbol.

\noindent Substituting this into the formula for \( C_f \), we obtain:
\[
C_f = \prod_{p} \left(1 - \frac{1 + \left( \frac{-163}{p} \right)}{p} \right) \left(1 - \frac{1}{p} \right)^{-1}.
\]

\noindent Computationally evaluating this product over all primes (up to a reasonably large cutoff), we find:
\[
\boxed{C_f \approx 3.3204}
\]

\noindent Using the Bateman-Horn conjecture, the expected number of prime values generated by the polynomial in the range \( n \in [0, 199] \) is:
\[
\pi_f(200) \approx C_f \cdot \int_2^{200} \frac{dt}{\log t} \approx C_f \cdot \frac{200}{\log 200} \approx 3.3204 \cdot \frac{200}{5.3} \approx \textbf{125.34}.
\]

\noindent However, the \textbf{actual number of primes} generated by the polynomial in this range is:
\[
\boxed{172}
\]

\vspace{0.5em}
\noindent This is significantly greater than the expected number, which empirically supports the fact that this polynomial is \textbf{prime-rich}, producing more primes than even the Bateman-Horn conjecture predicts on average. This exceptional behavior can be linked to the discriminant \( D = -163 \), which is one of the Heegner numbers associated with class number 1, historically known for their rich prime-generating behavior.

\vspace{1em}
\noindent\textbf{Example 2.} Consider the values \( Z = 2 \), \( H = 163 \), and \( k = 125 \), for \( n \in [0, 499] \). The corresponding polynomial becomes:
\[
f_{(Z,k,H)}(n) = n^2 - (2Zk - 1)n + \frac{(2Zk - 1)^2 + H}{4} = n^2 - 499n + 62291.
\]

\noindent This is a quadratic polynomial of the form \( f(n) = n^2 + an + b \), where:
\[
a = -499, \quad b = 62291.
\]

\noindent The discriminant is
\[
D = a^2 - 4b = (-499)^2 - 4 \cdot 62291 = 249001 - 249164 = -163.
\]

\noindent To compute the Bateman-Horn constant \( C_f \), we use the formula:
\[
C_f = \prod_{p} \left(1 - \frac{\omega_f(p)}{p} \right) \left(1 - \frac{1}{p} \right)^{-1},
\]
\noindent where \( \omega_f(p) \) is the number of solutions to \( f(n) \equiv 0 \mod p \). This depends on whether \( D = -163 \) is a quadratic residue modulo \( p \).

\noindent After computing the convergent infinite product up to a sufficiently large prime bound, we obtain the following:
\[
\boxed{C_f \approx 3.320376}
\]

\noindent To estimate the number of prime values of \( f(n) \) for \( n \in [0, 499] \), we use the Bateman-Horn heuristic:
\[
\pi_f(500) \sim \sum_{n=0}^{499} \frac{C_f}{\log f(n)}.
\]

\noindent Numerically evaluating this sum gives:
\[
\text{Expected number of primes: } \textbf{267.14}
\]

\noindent However, the actual number of primes for this polynomial in the range is:
\[
\boxed{368}
\]

\noindent This result again confirms the \textbf{exceptionally prime-rich} nature of this family of polynomials.

\vspace{1em}
\noindent\textbf{Example 3.} Now consider the values \( Z = 4 \), \( H = 163 \), and \( k = 125 \) with \( n \in [0, 999] \). The corresponding polynomial becomes:
\[
f_{(Z,k,H)}(n) = n^2 - (2Zk - 1)n + \frac{(2Zk - 1)^2 + H}{4} = n^2 - 999n + 249541.
\]

\noindent This is again a quadratic of the form \( f(n) = n^2 + an + b \), with:
\[
a = -999, \quad b = 249541.
\]

\noindent The discriminant of the polynomial is
\[
D = a^2 - 4b = (-999)^2 - 4 \cdot 249541 = 998001 - 998164 = -163.
\]

\noindent The discriminant is again negative, and specifically \( D = -163 \), which is the same as in the previous examples.

\noindent To compute the Bateman-Horn constant \( C_f \), we use the following.
\[
C_f = \prod_{p} \left(1 - \frac{\omega_f(p)}{p} \right) \left(1 - \frac{1}{p} \right)^{-1},
\]
where \( \omega_f(p) = 1 + \left( \frac{-163}{p} \right) \).

\noindent So,
\[
C_f = \prod_{p} \left(1 - \frac{1 + \left( \frac{-163}{p} \right)}{p} \right) \left(1 - \frac{1}{p} \right)^{-1}.
\]

\noindent Using a computational approximation up to a reasonable large prime bound, we get the following.
\[
 \boxed{C_f \approx 3.319940}
\]

\noindent The expected number of primes among \( f(n) \) for \( n \in [0, 999] \) is then approximated by:
\[
\text{Expected primes} \approx C_f \cdot \int_2^{1000} \frac{dt}{\log t} \approx 3.319940 \cdot \operatorname{Li}(1000)\] \[\therefore\text{Expected primes}\approx 3.319940 \cdot \frac{1000}{\log 1000} \approx \textbf{480.61}.
\]

\noindent However, the actual number of primes generated by this polynomial in the given range is:
\[
\boxed{ 652}
\]

\noindent The actual number of primes generated is significantly higher than the expected value under the Bateman-Horn conjecture, reinforcing the \textbf{prime-rich nature} of the polynomial \( f_{(Z,k,H)}(n) \).

\vspace{1em}

\noindent Similarly to examples \(1\), \(2\) and \(3\), several other polynomials of the form \( f_{(Z,k,H)}(n) \) were analyzed for their prime-generating behavior in defined ranges of \( n \). In each case, the actual number of primes produced by these polynomials significantly exceeded the expected count predicted by the Bateman-Horn conjecture. This consistent overperformance strongly suggests that all polynomials belonging to the family \( f_{(Z,k,H)}(n) \) exhibit a \textbf{prime-rich} property within suitable bounds of \( n \). Thus, this family of polynomials demonstrates exceptional potential in generating a large number of prime values and supports the hypothesis of their prime-generating efficiency.
The Bate-Horn conjecture gives an asymptotic estimate of the number of prime values that a polynomial (or a set of polynomials) can take \cite{batemanhorn}. A critical component of this estimate is the Bateman-Horn constant \( C_f \), defined as an infinite product on all primes. However, computing \( C_f \) directly becomes computationally intensive for large values of \( n \), as it requires evaluating \( \omega_f(p) \), the number of distinct roots modulo \( p \), for each prime \( p \). In the following section, we shall derive an approximate expression for \( C_f \) using the theory of quadratic residues, which will allow us to bypass evaluating each prime individually.

\section{Approximation of the Bateman-Horn Constant using Quadratic Residue Theory}

\noindent Since from eq.(\ref{4.2}), we know that the Bateman-Horn constant \( C_f \) for a polynomial \( f(n) \in \mathbb{Z}[n] \) is given by:
\[
C_f = \prod_{\substack{p\ \text{prime}}} \left( 1 - \frac{\omega_f(p)}{p} \right) \left(1 - \frac{1}{p}\right)^{-1},
\]
\noindent where \( \omega_f(p) \) is the number of distinct roots of \( f(n) \equiv 0 \pmod{p} \). For a quadratic polynomial, the value of \( \omega_f(p) \) depends on whether the discriminant is a quadratic residue modulo \( p \) \cite{quadraticresidues}.

\noindent For our polynomial
\[
f_{(Z,k,H)}(n) = n^2 - (2Zk - 1)n + \frac{(2Zk - 1)^2 + H}{4},
\]
The discriminant is \( D = (2Zk - 1)^2 - 4 \cdot \frac{(2Zk - 1)^2 + H}{4} = -H \). Hence, the number of modulo solutions \( p \) depends on whether \( -H \) is a quadratic residue mod \( p \).

\noindent We define

\begin{itemize}
    \item If \( -H \) is a \textbf{quadratic residue} \( \mod p \), then the polynomial has \textbf{2 distinct roots}, so \( \omega_f(p) = 2 \).
    \item If \( -H \) is a \textbf{quadratic non-residue} \( \mod p \), then the polynomial has \textbf{no root}, so \( \omega_f(p) = 0 \).
\end{itemize}

\noindent This allows us to divide the product for \( C_f \) as follows:
\[
C_f = \prod_{p \in \text{QR}} \left( \frac{p - 2}{p - 1} \right) \cdot \prod_{p \in \text{NQR}} \left( \frac{p}{p - 1} \right) \tag{4.1}\label{4.1}
\]

\noindent Taking natural logarithms on both sides:
\[
\log C_f = \sum_{p \in \text{QR}} \log\left( \frac{p - 2}{p - 1} \right) + \sum_{p \in \text{NQR}} \log\left( \frac{p}{p - 1} \right) \tag{4.2} \label{4.2}
\]

\noindent Practically, the sets of primes where \( -H \) is QR or NQR are randomly distributed \cite{quadraticresidues}. So, we introduce probabilities \( \Pr(\text{QR}) \) and \( \Pr(\text{NQR}) \), and write:
\[\boxed{
\log C_f = \Pr(\text{QR}) \cdot \sum_{p \in \text{QR}} \log\left( \frac{p - 2}{p - 1} \right) + \Pr(\text{NQR}) \cdot \sum_{p \in \text{NQR}} \log\left( \frac{p}{p - 1} \right) } \tag{4.3} \label{4.3}
\] 

\noindent Now, note the following approximations from Taylor series:
\[
\log\left( \frac{p - 2}{p - 1} \right) = \log\left( 1 - \frac{1}{p - 1} \right) \approx -\frac{1}{p - 1} + \varepsilon_1 \tag{4.4} \label{4.4}
\]
\[
\log\left( \frac{p}{p - 1} \right) = \log\left( 1 + \frac{1}{p - 1} \right) \approx \frac{1}{p - 1} + \varepsilon_2 \tag{4.5}\label{4.5}
\]
where \( \varepsilon_1 \) and \( \varepsilon_2 \) represent higher-order term errors.

\noindent Substituting eq.(\ref{4.4}) and eq.(\ref{4.5}) in eq.(\ref{4.3}) we get,
\[
\log C_f \approx -\Pr(\text{QR}) \cdot \sum_{p \in \text{QR}} \frac{1}{p - 1} + \Pr(\text{NQR}) \cdot \sum_{p \in \text{NQR}} \frac{1}{p - 1} + \epsilon \tag{4.6}\label{4.6} 
\]
where \(\epsilon = \varepsilon_1 +  \varepsilon_2 \)

\noindent Since \( p - 1 \approx p \) for large \( p \), we apply the approximation from Euler's product for the Riemann zeta function expressed in eq.(\ref{1.6}). Therfore, we obtain the following expression:

\[
\log C_f \approx -\Pr(\text{QR}) \cdot (\log \log x + \gamma) + \Pr(\text{NQR}) \cdot (\log \log x + \gamma),
+ \epsilon\]
\[ \therefore
\log C_f \approx (\Pr(\text{NQR}) - \Pr(\text{QR})) \cdot (\log \log x + \gamma') \tag{4.7}\label{4.7}
\]

\noindent Let \( \delta_P = |\Pr(\text{NQR}) - \Pr(\text{QR})| \), then:
\[
\log C_f \approx \delta_P \cdot (\log \log x + \gamma') \tag{4.8}\label{4.8}
\]

\noindent Through computational evaluation and fitting, it was found that \( \gamma' \approx 39.1751 \). Therefore, the final approximate expression becomes:
\[
\log C_f \approx \delta_P \cdot (\log \log x + 39.1751) 
\]
or equivalently,
\[
C_f \approx \exp\left( \delta_P \cdot (\log \log x + 39.1751) \right) \tag{4.9}\label{4.9}
\]

\noindent [\textbf{Note: It was empirically observed that for the family \( f_{Z,k,H}(n) \), \( \delta_P \) tends to stabilize around a value of \( 0.03023 \).}]\cite{quadraticresidues}\\
After substituting \( \delta_P \) =  \( 0.03023 \) in eq.(\ref{4.9}), the final approximation becomes the following:
\[\boxed{
C_f \approx \exp\left( 0.03023 \cdot (\log \log x + 39.1751) \right)} \tag{4.10}\label{4.10}
\]

\noindent This approximation tells us that the Bateman-Horn constant \( C_f \), which governs the density of primes generated by our polynomial, can be efficiently estimated without computing products over all primes. The simplification uses the probability of \( -H \) being a quadratic residue or not, and relates it directly to known logarithmic growth rates and constants. This is particularly useful for large \( n \)-ranges, where full prime factor analysis becomes computationally heavy. The approximation retains the essential characteristics of \( C_f \), giving us a fast and insightful tool to evaluate the prime-generating nature of our polynomial.

Let us consider a few examples to gain a clearer understanding of the approximation of \( C_f \) and the precision of this estimation.\\

\noindent \textbf{Example 1:} Consider the prime-rich polynomial from the family \( f_{Z,k,H}(n) \), with parameters: \( Z = 1 \),  \( k = 80 \), \( H = 163 \), \( n \in [0,99] \).\\

\noindent The polynomial is defined as:
\[ f_{(Z,k,H)}(n) = n^2 - An + C \quad \text{where} \quad A = 2Zk - 1 = 159, \quad C = \frac{A^2 + H}{4} = 6361\]
\[\therefore f_{(1,80,163)}(n) = n^2 - 159n + 6361 \]

\noindent We compute the actual Bateman-Horn constant over the first 100 values (\( n = 0 \) to \( 99 \)):
\[ C_f = \prod_{p\text{ prime}} \left(1 - \frac{\omega_f(p)}{p}\right) \left(1 - \frac{1}{p}\right)^{-1} \]
where \( \omega_f(p) \) is the number of solutions to \( f(n) \equiv 0 \mod p \).

\noindent Using computational evaluation, the actual value is found to be:
\[ \boxed{C_f = 3.3259186569785313} \]

\noindent The discriminant of the polynomial is:
\[ D = A^2 - 4C = 159^2 - 4 \cdot 6361 = -163 \]

\noindent We now compute the number of primes \( p \leq 6361 \) (i.e., maximum prime factor of \( f(n) \)) for which \( \left( \frac{D}{p} \right) = 1 \) (QR) or \( -1 \) (NQR). We find the following.
\[ \delta_P = \frac{|\#QR - \#NQR|}{\#QR + \#NQR} \approx 0.03023 \]

\noindent The approximate value of \( C_f \) will be:
\[ C_f \approx \exp\left( \delta_P \cdot \left( \log(\log x) + 39.1751 \right) \right) \]
where \( x = 6361 \) is the maximum prime divisor of any \( f(n) \).

\[
\therefore
C_f \approx \exp(0.03023 \cdot (\log(\log(6361)) + 39.1751)) \]
\[C_f \approx  \exp(0.03023 \cdot 39.75528632) \approx \boxed{3.326106183}
\]

\noindent The approximated value \( C_f \approx 3.326106183 \) is extremely close to the actual constant \( C_f = 3.3259186569785313 \), thus validating the accuracy of the \( \delta_P \)-based heuristic for this polynomial.\\

\noindent \textbf{Example 2:} Consider the prime-rich polynomial from the family \( f_{Z,k,H}(n) \), with parameters: \( Z = 2 \),  \( k = 50 \), \( H = 163 \), \( n \in [0,200] \).\\

\noindent The polynomial is defined as:
\[ f_{(Z,k,H)}(n) = n^2 - An + C \quad \text{where} \quad A = 2Zk - 1 = 199, \quad C = \frac{A^2 + H}{4} = 9941 \]
\[\therefore f_{(2,50,163)}(n) = n^2 - 199n + 9941 \]

\noindent We compute the actual Bateman-Horn constant over the first 201 values (\( n = 0 \) to \( 200 \)):
\[ C_f = \prod_{p\text{ prime}} \left(1 - \frac{\omega_f(p)}{p}\right) \left(1 - \frac{1}{p}\right)^{-1} \]
where \( \omega_f(p) \) is the number of solutions to \( f(n) \equiv 0 \mod p \).

\noindent Using computational evaluation, the actual value is found to be:
\[ \boxed{C_f = 3.3283621082005626} \]

\noindent The discriminant of the polynomial is:
\[ D = A^2 - 4C = 199^2 - 4 \cdot 9941 = -163 \]

\noindent We now compute the number of primes \( p \leq 9941 \) (i.e., maximum prime factor of \( f(n) \)) for which \( \left( \frac{D}{p} \right) = 1 \) (QR) or \( -1 \) (NQR). We find the following.
\[ \delta_P = \frac{|\#QR - \#NQR|}{\#QR + \#NQR} \approx 0.03023 \]

\noindent The approximate value of \( C_f \) will be:
\[ C_f \approx \exp\left( \delta_P \cdot \left( \log(\log x) + 39.1751 \right) \right) \]
where \( x = 9941 \) is the maximum prime divisor of any \( f(n) \).

\[
\therefore
C_f \approx \exp(0.03023 \cdot (\log(\log(9941)) + 39.1751)) \]
\[C_f \approx  \exp(0.03023 \cdot 39.77688088) \approx \boxed{3.328278186}
\]

\noindent The approximated value \( C_f \approx 3.328278186 \) is extremely close to the actual constant \( C_f = 3.3283621082005626 \), thus validating the accuracy of the \( \delta_P \)-based heuristic for this polynomial.\\

\noindent \textbf{Example 3:} Consider the prime-rich polynomial from the family \( f_{Z,k,H}(n) \), with parameters: \( Z = 1 \),  \( k = 160 \), \( H = 163 \), \( n \in [0,200] \).\\

\noindent The polynomial is defined as:
\[ f_{(Z,k,H)}(n) = n^2 - An + C \quad \text{where} \quad A = 2Zk - 1 = 319, \quad C = \frac{A^2 + H}{4} = 25481 \]
\[\therefore f_{(1,160,163)}(n) = n^2 - 319n + 25481 \]

\noindent We compute the actual Bateman-Horn constant over the first 201 values (\( n = 0 \) to \( 200 \)):
\[ C_f = \prod_{p\text{ prime}} \left(1 - \frac{\omega_f(p)}{p}\right) \left(1 - \frac{1}{p}\right)^{-1} \]
where \( \omega_f(p) \) is the number of solutions to \( f(n) \equiv 0 \mod p \).

\noindent Using computational evaluation, the actual value is found to be:
\[ \boxed{C_f = 3.3206295784604256} \]

\noindent The discriminant of the polynomial is:
\[ D = A^2 - 4C = 319^2 - 4 \cdot 25481 = -163 \]

\noindent We now compute the number of primes \( p \leq 25481 \) (i.e., maximum prime factor of \( f(n) \)) for which \( \left( \frac{D}{p} \right) = 1 \) (QR) or \( -1 \) (NQR). We find the following.
\[ \delta_P = \frac{|\#QR - \#NQR|}{\#QR + \#NQR} \approx 0.03023 \]

\noindent The approximate value of \( C_f \) will be:
\[ C_f \approx \exp\left( \delta_P \cdot \left( \log(\log x) + 39.1751 \right) \right) \]
where \( x = 25471 \) is the maximum prime divisor of any \( f(n) \).

\[
\therefore
C_f \approx \exp(0.03023 \cdot (\log(\log(25471)) + 39.1751)) \]
\[C_f \approx  \exp(0.03023 \cdot 39.81914903) \approx \boxed{3.332533665}
\]

\noindent The approximated value \( C_f \approx 3.332533665 \) is extremely close to the actual constant \( C_f = 3.3206295784604256 \), thus validating the accuracy of the \( \delta_P \)-based heuristic for this polynomial.\\

\paragraph{}
After defining an approximate value for the Bateman-Horn constant \( C_f \), let us now dive into an approximation that shall lead to determining the maximum number of prime-generating polynomials for a given range of \( n \). This involves deriving an optimal value for \( Zk \) such that the corresponding polynomial exhibits the highest prime density.

\section{Derivation of Optimal \texorpdfstring{$Zk$}{Zk} for Maximum Prime Generation}

\noindent As the prime-rich polynomial family expressed in eq.(\ref{1.2}) is a quadratic polynomial, and the vertex of a general quadratic of the form \( an^2 + bn + c \) is given by:
\[
n = -\frac{b}{2a}
\]
\noindent Comparing coefficients, we get:
\[
a = 1, \quad b = -(2Zk - 1), \quad c = \frac{(2Zk - 1)^2 + H}{4}
\]
\noindent Therefore, the vertex of \( f_{(Z,k,H)}(n) \) occurs at:
\[
n = -\frac{- (2Zk - 1)}{2} = \frac{2Zk - 1}{2} = K(assume) \tag{5.1} \label{5.1}
\]
\noindent Since \( K \in \mathbb{N} \), it is the point of symmetry of the polynomial. We now reparametrize the polynomial in terms of \( K \), rather than \( Z \) and \( k \).

\noindent From eq.(\ref{5.1}), we have:
\[
f_K(n) = n^2 - 2Kn + K^2 + \frac{H}{4}
\]
\noindent This can be rewritten as:
\[
f_K(n) = (n - K)^2 + \frac{H}{4} \tag{5.2}\label{5.2}
\]
\noindent Let \( c = \frac{H}{4} \), a constant (since \( H \) is a fixed Heegner number), so that:
\[
f_K(n) = (n - K)^2 + c \tag{5.3}\label{5.3}
\]

\noindent According to the modified Bateman-Horn conjecture, the expected number of primes generated by the polynomial \( f_K(n) \) over the interval \( [n_1, n_2] \) is asymptotically given by:
\[
\pi_{\text{expected}}(K) \propto \int_{n_1}^{n_2} \frac{1}{\log(f_K(n))} \, dn
\]
\noindent Substituting eq.(\ref{5.3}), we get:
\[
\pi_{\text{expected}}(K) \propto \int_{n_1}^{n_2} \frac{1}{\log((n - K)^2 + c)} \, dn
\]

\[
\therefore \pi_{\text{expected}}(K) = C_f(K) \cdot \int_{n_1}^{n_2} \frac{1}{\log((n - K)^2 + c)} \, dn \tag{5.4}\label{5.4}
\]
where, \( C_f(K) \) is the Bateman-Horn constant for \( f_K(n) \).
\noindent To find the value of \( K \) that maximizes \( \pi_{\text{expected}}(K) \), we differentiate with respect to \( K \) using Chain Rule:
\[\
\frac{d}{dK} \pi_{\text{expected}}(K) = C_f'(K) \cdot \int_{n_1}^{n_2} \frac{1}{\log((n - K)^2 + c)} \, dn + C_f(K) \cdot \frac{d}{dK} \left( \int_{n_1}^{n_2} \frac{1}{\log((n - K)^2 + c)} \, dn \right)
\]
\noindent since, \(C_f(K) =\) Bateman-Horn constant,
\(\therefore C_f'(K) = 0\)
\[ \therefore \frac{d}{dK} \pi_{\text{expected}}(K) = C_f(K) \cdot \frac{d}{dK} \left( \int_{n_1}^{n_2} \frac{1}{\log((n - K)^2 + c)} \, dn \right) \tag{5.5}\label{5.5}
\]

\noindent Using \textit{Leibniz’s rule} \cite{Leibnitz} for differentiation under the integral sign:
\[
\frac{d}{dK} \left( \int_{n_1}^{n_2} \frac{1}{\log((n - K)^2 + c)} \, dn \right) = \int_{n_1}^{n_2} \frac{d}{dK} \left( \frac{1}{\log((n - K)^2 + c)} \right) \, dn
\]

\noindent Compute the inner derivative:
\[
\frac{d}{dK} \left( \frac{1}{\log((n - K)^2 + c)} \right) = \frac{-2(n - K)}{((n - K)^2 + c) \cdot (\log((n - K)^2 + c))^2}
\]

\noindent Let us use substitution here:
\[
u = n - K \quad \Rightarrow \quad n = u + K \quad \text{and} \quad dn = du
\]
\noindent Then the limits become:
\[
n = n_1 \Rightarrow u = n_1 - K, \quad n = n_2 \Rightarrow u = n_2 - K
\]

\noindent The integral becomes:
\[
\int_{n_1}^{n_2} \frac{1}{\log((n - K)^2 + c)} \, dn = \int_{n_1 - K}^{n_2 - K} \frac{1}{\log(u^2 + c)} \, du
\]

\noindent Let
\[
F(K) = \int_{n_1 - K}^{n_2 - K} \frac{1}{\log(u^2 + c)} \, du
\]

\noindent Using Leibniz's rule:
\[
\frac{dF}{dK} = \frac{d}{dK} \left( \int_{a(K)}^{b(K)} f(u) \, du \right) = f(b(K)) \cdot \frac{db}{dK} - f(a(K)) \cdot \frac{da}{dK}
\]

\noindent Where:
\begin{align*}
a(K) &= n_1 - K \quad \Rightarrow \quad \frac{da}{dK} = -1 \\
b(K) &= n_2 - K \quad \Rightarrow \quad \frac{db}{dK} = -1 \\
f(u) &= \frac{1}{\log(u^2 + c)}
\end{align*}

\noindent Hence,
\[
\frac{dF}{dK} = -\frac{1}{\log((n_2 - K)^2 + c)} + \frac{1}{\log((n_1 - K)^2 + c)}
\]

\[
\frac{d}{dK} \left( \int_{n_1}^{n_2} \frac{1}{\log((n - K)^2 + c)} \, dn \right) = \frac{1}{\log((n_1 - K)^2 + c)} - \frac{1}{\log((n_2 - K)^2 + c)}
\]

\[\frac{d}{dK} \pi_{\text{expected}}(K) = C_f(K) \left[ \frac{1}{\log((K - n_1)^2 + c)} - \frac{1}{\log((K - n_2)^2 + c)} \right] \tag{5.6}\label{5.6}
\]

\noindent To find the extremum (maximum), set eq.(\ref{5.6}) to zero:
\[
\frac{1}{\log((K - n_1)^2 + c)} = \frac{1}{\log((K - n_2)^2 + c)}
\]
\[
\therefore\log((K - n_1)^2 + c) = \log((K - n_2)^2 + c)
\]
\noindent Taking exponentials on both sides:
\[
(K - n_1)^2 + c = (K - n_2)^2 + c\]
\[(K - n_1)^2 = (K - n_2)^2
\]
\[
(K^{2} - 2.K.n_{1} + n_{1}^{2}) = (K^{2} - 2.K.n_{2} + n_{2}^{2})
\]

\[2K = n_1 + n_2 \]

\[\therefore
K = \frac{n_1 + n_2}{2} \tag{5.7}\label{5.7}
\]

\noindent Substituting, \( K = \frac{2Zk - 1}{2} \) in eq.(\ref{5.7}) we obtain:
\[
\frac{2Zk - 1}{2} = \frac{n_1 + n_2}{2} \Rightarrow 2Zk - 1 = n_1 + n_2
\Rightarrow Zk = \frac{n_1 + n_2 + 1}{2}
\]

\noindent Thus, the approximate value of \( Zk \) that maximizes the number of primes generated by the polynomial over the interval \( [n_1, n_2] \) is:
\[\boxed{
Zk \approx \frac{n_1 + n_2 + 1}{2}} 
\]

\noindent However, 2 cases arise during the approximation of \(Z.k\). They are expressed as follows:

\noindent
\textbf{Case 1:} \( Zk \) is an integer if \( (n_1 + n_2 + 1) \) is even. In this case, the approximation directly yields a valid integer value for \( Zk \), which can be used in the polynomial without any modification.

\noindent
\textbf{Case 2:} \( Zk \) is not an integer if \( (n_1 + n_2 + 1) \) is odd. Since \( Zk \) must be a positive integer for the polynomial to retain its prime-rich structure, we take the floor or the rounded value of \( \frac{n_1 + n_2 + 1}{2} \) to obtain an integral value of \( Zk \) that can be practically used in computation.

\noindent Hence, the general approximation can be expressed as:

\[\boxed{
Zk \approx \left\lfloor \frac{n_1 + n_2 + 1}{2} \right\rfloor
\quad \text{or} \quad Zk \approx \operatorname{round}\left( \frac{n_1 + n_2 + 1}{2} \right)} \tag{5.8} \label{5.8}\]

Eqn.(\ref{5.8}) provides an estimate for selecting a \( Zk \) such that the corresponding polynomial \( f_{(Z,k,H)}(n) \) produces the maximal number of primes in the given range. Note that this is an \textit{approximation}, as the Bateman-Horn conjecture only predicts \textit{asymptotic behavior}. For larger ranges of \( n \), the actual number of primes may differ from the expected due to local fluctuations.
Let us have a clear understanding of this approximation with the help of a few examples.\\

\noindent \textbf{Example 1:} Let the given range for \( n \) be \( [0, 99] \). Here, \( n_1 = 0 \) and \( n_2 = 99 \). Therefore from eq.(\ref{5.8}),
\[\boxed{
Zk \approx \left\lfloor \frac{n_1 + n_2 + 1}{2} \right\rfloor = \left\lfloor \frac{0 + 99 + 1}{2} \right\rfloor = 50}
\]

\noindent Substituting \( Zk = 50 \) into the polynomial
\[
f_{(Z,k,H)}(n) = n^2 - (2Zk - 1)n + \frac{(2Zk - 1)^2 + H}{4}
\]
\noindent with \( H = 163 \) (a Heegner number), we get:
\[
f_{(Z,k,H)}(n) = n^2 - 99n + 2491
\]

\noindent On computational analysis of the total number of primes generated by this polynomial for \( n \in [0, 99] \), we find that it generates \textbf{92 primes and 8 non-primes}, which is very close to the maximum number of primes (\textbf{95}) achievable in this range. This validates the accuracy of the approximation.\\

\noindent \textbf{Example 2:} Let the given range for \( n \) be \( [21, 95] \). Here, \( n_1 = 21 \) and \( n_2 = 95 \). Therefore from eq.(\ref{5.8}),
\[ \boxed{
Zk \approx \left\lfloor \frac{21 + 95 + 1}{2} \right\rfloor = \left\lfloor \frac{117}{2} \right\rfloor = 58} \quad
 \text{or} \quad
\boxed{Zk \approx \text{round}\left( \frac{21 + 95 + 1}{2} \right) = 59}
\]

\noindent\textbf{Case 1: \( Zk = 58 \)}\\
\noindent Substituting into the polynomial:
\[
f_{(Z,k,H)}(n) = n^2 - 115n + 3347
\]

\noindent On computational analysis of this polynomial over the interval \( [21, 95] \), we find that it generates \textbf{75 primes and 0 non-primes}, which is the maximum possible in this range.

\noindent \textbf{Case 2: \( Zk = 59 \)}\\
\noindent Substituting into the polynomial:
\[
f_{(Z,k,H)}(n) = n^2 - 117n + 3463
\]

\noindent Again, over the interval \( [21, 95] \), the polynomial generates \textbf{75 primes and 0 non-primes}, achieving the same maximum prime count.

\noindent These examples confirm the approximation,
\[
Zk \approx \left\lfloor \frac{n_1 + n_2 + 1}{2} \right\rfloor \quad \text{or} \quad Zk \approx \text{round}\left( \frac{n_2 + n_1 + 1}{2} \right)
\]
\noindent is a valid and accurate method to estimate the value of \( Zk \) that maximizes the number of primes generated by the polynomial \( f_{(Z,k,H)}(n) \) over a specified range of \( n \).

After understanding the approximated most favorable values of $Z \cdot k$ for any given range of $n \in [n_1, n_2]$, 
we now arrive at a very interesting and significant portion of our study. 
What makes our polynomial $f_{(Z,k,H)}(n)$ particularly special is its profound correlation with the celebrated 
\textbf{Riemann Zeta function} and how its structure closely aligns with the central claim of the \textbf{Riemann Hypothesis}. 
In what follows, we delve deeply into this connection, interpret the roots of the polynomial in the complex plane, 
and reveal how they mirror the non-trivial zeros of $\zeta(s)$ in form and symmetry.

\section{Correlation Between the Prime-Rich Polynomial \texorpdfstring{$f_{(Z,k,H)}(n)$}{f(Z,k,H)(n)} and the Riemann Hypothesis}

\noindent As mentioned in eq.(\ref{1.5}), the \textbf{Riemann zeta function} is defined for complex numbers $s = \alpha + i\beta$ with $\alpha > 1$ as:
\[
\zeta(s) = \sum_{n=1}^{\infty} \frac{1}{n^s}
\]
\noindent It can be continued analytically with the rest of the complex plane, except for a simple pole at $s = 1$.

\noindent Similarly, the \textbf{Riemann Hypothesis} postulates that all non-trivial zeros of the Riemann zeta function lie on \emph{critical line}:
\[
\text{Re}(s) = \frac{1}{2}
\]
\noindent That is, if $\zeta(s) = 0$ and $s$ is a non-trivial zero, then $s = \frac{1}{2} + i\beta$ for some real number $\beta$.
This hypothesis has profound implications for the distribution of prime numbers. The truth of the hypothesis would confirm that primes are distributed as regularly as possible, in accordance with the predictions of the Prime Number Theorem \cite{primenumbertheorem}.

\noindent After a thorough analysis we discovered that the polynomials constructed from eq.(\ref{1.2}) are \textit{prime-rich}, i.e., for suitable choices of parameters, they generate large sequences of prime numbers for consecutive integer values of $n$.

\noindent The roots of the polynomial $f_{(Z,k,H)}(n)$ are given by solving the quadratic:
\[
n = \frac{(2Zk - 1) \pm \sqrt{(-(2Zk - 1))^{2} -4.[\frac{(2Zk - 1)^{2}+H}{4}]}}{2}
\]
\[
\therefore n = \frac{(2Zk - 1) \pm \sqrt{-H}}{2}
\]
\noindent This yields complex roots of the form:
\[
n = \frac{2Zk - 1}{2} \pm i\frac{\sqrt{H}}{2}
\]

\noindent We can factorize the polynomial as follows:
\[\boxed{
f_{(Z,k,H)}(n) = \left(n - \left(\frac{2Zk - 1}{2} + i\frac{\sqrt{H}}{2} \right) \right)\left(n - \left(\frac{2Zk - 1}{2} - i\frac{\sqrt{H}}{2} \right) \right)} \tag{6.1} \label{6.1}
\]

\noindent Eq.(\ref{6.1}) clearly resembles the general form of a complex number expressed in terms of its real and imaginary parts expressed as:
\[
n = \alpha \pm i\beta
\]
\noindent which mirrors the structure of the non-trivial zeros of the Riemann zeta function $s = \alpha + i\beta$.

\noindent If we set $Z = k = 1$, then the polynomial simplifies to:
\[
f_{(1,1,H)}(n) = n^2 - n + \frac{1 + H}{4}
\]
\noindent and the roots become:
\[ \boxed{
n = \frac{1}{2} \pm i\frac{\sqrt{H}}{2}} \tag{6.2} \label{6.2}
\]

\noindent In eq.(\ref{6.2}), the real part of the roots is exactly $\frac{1}{2}$ which resembles the critical line $\text{Re}(s) = \frac{1}{2}$ proposed in the Riemann Hypothesis. This is a profound observation. It suggests a structural similarity between the roots of the prime-rich polynomial and the non-trivial zeros of $\zeta(s)$.
This connection is not merely coincidental. It suggests that the algebraic structure underlying certain families of prime-generating polynomials may be inherently tied to the distribution of prime numbers which becomes the central concern of the Riemann Hypothesis.

\begin{itemize}
    \item The imaginary part $\sqrt{H}/2$ plays a role analogous to the imaginary part $\beta$ in the zeta zeros.
    \item The real part $\alpha = 1/2$ is a shared symmetry in both systems.
    \item The factorization of the polynomial resembles the structure of the product in the Euler product representation of $\zeta(s)$:
    \[
    \zeta(s) = \prod_{p\ \text{prime}} \left(1 - \frac{1}{p^s} \right)^{-1}
    \]
    where the primes are tied to the zeros of $\zeta(s)$ through their logarithmic distribution.
\end{itemize}

\noindent Thus, the structure of $f_{(Z,k,H)}(n)$ may capture deep analytical properties of prime distributions. Although this does not prove the Riemann Hypothesis, it strongly suggests a harmonious algebraic analogy, and possibly a new heuristic framework, for understanding the location of zeta zeros via prime-rich polynomial structures.

\noindent The family of polynomials $f_{(Z,k,H)}(n)$ not only generates primes but also encodes complex root structures that align with the critical features of the Riemann zeta function. This correspondence between the roots of polynomials and the non-trivial zeros of $\zeta(s)$ adds rich mathematical depth and may inspire further theoretical development bridging algebraic and analytic number theory. It reaffirms the remarkable unity in the world of prime numbers, whether approached through elementary polynomials or deep complex-analytic functions.\\

Until now, we have focused more on the theoretical concepts surrounding our polynomial $f_{(Z,k,H)}(n)$, 
particularly emphasizing its prime-rich property in a structured way through various comparisons, theorems, and heuristic conjectures. 
We have demonstrated its strong potential to generate prime numbers effectively and densely in selected ranges of $n$. 
Now, we shall broaden our perspective by delving deep into the diverse areas where this polynomial can prove to be a turning point in the field of number theory and beyond.

Prime numbers play a pivotal role in real-world applications, most notably in cryptography, digital security systems, 
hashing algorithms, random number generation, error-correcting codes, and blockchain technology. 
However, the unpredictability and irregular distribution of primes have historically posed significant challenges in practical implementations. Efficient prime generation and prediction have always been bottlenecks in secure communication systems and data encryption algorithms.

This is where the polynomial $f_{(Z,k,H)}(n)$ becomes a game changer. Its structured form and tunable parameters allow for customized generation of prime-dense sequences, addressing the hurdles of randomness and computational inefficiency. By encapsulating deep mathematical symmetry and offering a unifying framework, this polynomial has the potential to revolutionize the way primes are generated, verified, and applied in real-life systems.

With this foundational understanding, we now proceed to explore the real-world applications of the polynomial $f_{(Z,k,H)}(n)$ and how it can be used to solve long-standing challenges in mathematics, computer science, and engineering.

\section{The Role of the Polynomial \texorpdfstring{$f_{(Z,k,H)}(n)$}{f(Z,k,H)(n)} in Enhancing Modern Cryptographic Frameworks}

In traditional asymmetric encryption schemes such as RSA, the core of security lies in the computational difficulty of factoring the product of two large prime numbers. However, as computational power advances and especially with the looming development of quantum algorithms such as the Shor’s algorithm, there is an increasing need to reinforce the generation and structure of these primes to maintain cryptographic integrity. In this context, the polynomial $f_{(Z,k,H)}(n)$ offers a powerful and structured approach to the generation of prime numbers and enhanced key protection, making it a novel contribution to cryptographic systems.

Unlike standard methods in which two large primes are randomly selected, our method begins by selecting these primes through a defined, secure mathematical process. First, two random but large and distinct values for $Z$ and $k$ are selected. These values are then used to evaluate the expression:

\[
\frac{(2Zk - 1)^2 + H}{4} \tag{7.1} \label{7.1}
\]

\noindent We check whether $-H$ is a quadratic residue modulo eq.(\ref{8.1}) using the Jacobi symbol, denoted as $\left( \frac{-H}{m} \right)$, where $m = \frac{(2Zk - 1)^2 + H}{4}$. If the Jacobi test yields $+1$, we continue to perform a single iteration of the Miller-Rabin primality test on $m$ to confirm its primality. If $m$ fails either test, we discard the current pair $(Z,k)$1 and select new high values, repeating the process. Since the polynomial $f_{(Z,k,H)}(n)$ is inherently prime rich, the probability of finding suitable pairs $(Z,k)$ is high, ensuring efficient selection without compromising security.

\noindent This process is performed twice to obtain two distinct prime numbers using two pairs: $(Z_1, k_1)$ and $(Z_2, k_2)$. Now, observe the nature of the polynomial:

\[
f_{(Z,k,H)}(n) = n^2 - (2Zk - 1)n + \frac{(2Zk - 1)^2 + H}{4}
\]

\noindent At $n = 0$, we get:

\[ \boxed{
f_{(Z,k,H)}(0) = \frac{(2Zk - 1)^2 + H}{4}} \tag{7.2}\label{7.2}
\]

\noindent Therefore, it is not necessary to evaluate the polynomial for multiple $n$ values; evaluating it at $n = 0$ is sufficient. Thus, the two primes used to generate the public key are:

\[\boxed{
p_1 = f_{(Z_1,k_1,H)}(0), \quad p_2 = f_{(Z_2,k_2,H)}(0)} \tag{7.3}\label{7.3}
\]

\noindent These are then multiplied to produce the modulus $N = p_1 \cdot p_2$, which is published as part of the public key.

\noindent Now, suppose that an attacker successfully factors $N$ into $p_1$ and $p_2$. The attacker might attempt to deduce the private keys by solving the reverse of the expression:

\[\boxed{
p = \frac{(2Zk - 1)^2 + H}{4}
\Rightarrow (2Zk - 1)^2 = 4p - H
\Rightarrow Zk = \frac{\sqrt{4p - H} + 1}{2}} \tag{7.4} \label{7.4}
\]

\noindent However, eq.(\ref{7.4}) only produces the product $Zk$, not the individual values of $Z$ and $k$. As $Z$ and $k$ are distinct and large, recovering both from their product involves computing all possible factor pairs of $Zk$, an infeasible task due to the high number of permutations and combinations involved. Even if the attacker successfully finds one valid pair, they would still need to find the second pair, repeating the entire effort and effectively doubling the complexity.

To offer a real-world analogy: imagine a digital vault locked with two complex padlocks. These padlocks are not ordinary, they are designed based on intricate mechanisms derived from a mathematical blueprint. Even if someone breaks the outermost security to see what the locks look like (i.e., factors the primes), they still have to figure out the inner gear mechanisms (the exact $(Z,k)$ pairs) to open it. Without precise knowledge of how each padlock is internally structured, any attempt to open them would require testing thousands, if not millions, of internal configurations. The odds are daunting.

What truly differentiates this encryption approach is that the private key is not merely the primes, but the very structure, the $(Z,k)$ pairs, that generated those primes. This makes the system inherently more secure, even in a hypothetical future where prime factorization becomes trivial. In such a future, our method ensures that, without complete knowledge of the originating polynomial structure, unauthorized decryption remains computationally impractical. Thus, the polynomial $f_{(Z,k,H)}(n)$ not only enhances security today but also anticipates the cryptographic challenges of tomorrow with a mathematically fortified framework.

To demonstrate the practical implementation of this technique, the following algorithm outlines the structured generation of private and public keys using the polynomial $f_{(Z,k,H)}(n)$. This not only ensures the cryptographic strength of the keys, but also accelerates the generation of large primes with embedded structure for enhanced security.
This algorithm generates strong and structured prime numbers using the prime-rich polynomial $f_{(Z,k,H)}(n) = \frac{(2Zk - 1)^2 + H}{4}$, where $H$ is a hardcoded discriminant. The core function \texttt{GenerateFastStructuredPrime} randomly selects large integers $Z$ and $k$ within a predefined range, computes an intermediate base value $bk = 2Zk - 1$, and evaluates the candidate prime number using the polynomial. It then verifies the primality of the result by ensuring it meets a specified bit-length threshold and passes both the Jacobi symbol test $J(-H, p) = 1$ and the Miller-Rabin probabilistic primality test. If valid, the tuple $(Z, k, bk, p)$ is returned. This process is repeated to obtain two distinct primes $(p_1, p_2)$ along with their corresponding private parameters, ensuring uniqueness and non-null values. Finally, the public key is computed as the product $P = p_1 \cdot p_2$. By combining number-theoretic rigor and cryptographic strength, this algorithm ensures the generation of large, structured, and secure prime numbers.

\begin{algorithm}[H]
\caption{Structured Prime Generation and Key Construction using the Prime-Rich Polynomial \texorpdfstring{$f_{(Z,k,H)}(n)$}{f_{(Z,k,H)}(n)}}
\begin{algorithmic}[1] 
\State \textbf{Input:} Hardcoded discriminant $H = 163$, bit size requirement $bits = 200$, maximum attempts $max\_attempts = 10000$
\State \textbf{Output:} Two strong private keys $(Z_1, k_1)$ and $(Z_2, k_2)$, their primes $p_1, p_2$, and the public key $P = p_1 \cdot p_2$

\Function{GenerateFastStructuredPrime}{$H, bits, max\_attempts$}
    \For{$attempt = 1$ to $max\_attempts$}
        \State Randomly choose $Z \in [10^{80}, 10^{85}]$
        \State Randomly choose $k \in [10^{80}, 10^{85}]$
        \State Compute $bk \gets 2Zk - 1$
        \State Compute candidate prime: $p \gets \frac{bk^2 + H}{4}$
        \If{$\text{bit\_length}(p) \geq bits$}
            \If{Jacobi Symbol $J(-H, p) = 1$ \textbf{and} $p$ passes Miller-Rabin test}
                \State \Return $(Z, k, bk, p)$
            \EndIf
        \EndIf
    \EndFor
    \State \Return \textbf{Failure}
\EndFunction

\State Start timer
\State $(Z_1, k_1, bk_1, p_1) \gets$ \Call{GenerateFastStructuredPrime}{$H, bits, max\_attempts$}
\State $(Z_2, k_2, bk_2, p_2) \gets$ \Call{GenerateFastStructuredPrime}{$H, bits, max\_attempts$}

\While{$p_1 = p_2$ \textbf{or} any of $Z_1, k_1, p_1, Z_2, k_2, p_2$ is \textbf{None}}
    \State $(Z_2, k_2, bk_2, p_2) \gets$ \Call{GenerateFastStructuredPrime}{$H, bits, max\_attempts$}
\EndWhile

\State Compute public key: $P \gets p_1 \cdot p_2$
\State Stop timer and calculate execution time
\State \textbf{Output:} $(Z_1, k_1, p_1), (Z_2, k_2, p_2), P$ 
\end{algorithmic}

\end{algorithm}

Fig.(\ref{f6}) provides a visual representation of the structured prime generation algorithm used in the asymmetric encryption process. It clearly outlines the iterative process of generating two large primes \( p_1 \) and \( p_2 \) using the polynomial \( f_{(Z,k,H)}(n) = \frac{(2Zk - 1)^2 + H}{4} \), ensuring that each candidate satisfies specific primality conditions, such as the Jacobi symbol and Miller–Rabin test. The diagram also highlights the rejection and regeneration loop, which reinforce the robustness and efficiency of the algorithm. Ultimately, the flow chart emphasizes how the deterministic structure and randomness together contribute to generating strong cryptographic keys.

\begin{figure}[H]
    \centering
    \includegraphics[width=0.5\linewidth]{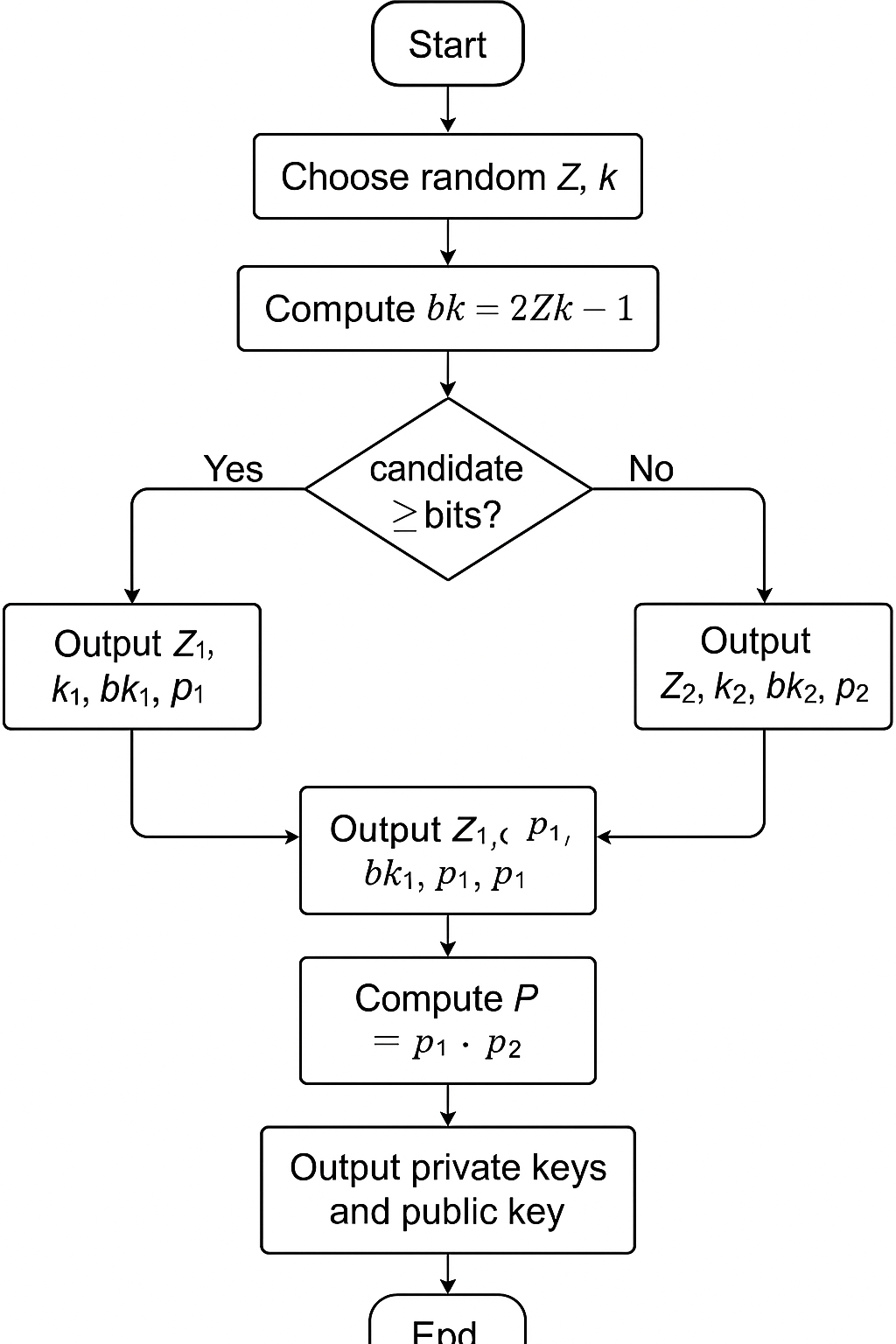}
    \caption{Flow chart of Algorithm 1}
    \label{f6}
\end{figure}

To evaluate the efficiency and security offered by the enhanced encryption method using polynomial \( f_{(Z,k,H)}(n) \), a comparative analysis against the widely used RSA algorithm was carried out. The following table summarizes the observed differences based on key metrics and execution performance.\\

\noindent The comparative analysis between the Enhanced Asymmetric Encryption and the Traditional RSA Algorithm shown in fig.(\ref{table3}), highlights several key distinctions in terms of key structure, bit length, and performance. The Enhanced Asymmetric Encryption method utilizes a more complex private key structure, involving two sets of private keys: \(Z_1, k_1\) and \(Z_2, k_2\), whereas the Traditional RSA Algorithm does not disclose such internal private key parameters. The enhanced method also features longer and more diversified prime numbers \(P_1\) and \(P_2\), with larger numerical values compared to the RSA primes, indicating a higher level of cryptographic complexity. Additionally, the Enhanced scheme uses Jacobi symbol checks for both \(P_1\) and \(P_2\), which are not applicable in the RSA setup.
In terms of key lengths, the Enhanced method yields a public key bit length of \textbf{2253 bits}, slightly greater than RSA's \textbf{2046 bits}. The private key lengths also differ significantly: for the enhanced system, the bit lengths for \(Z_1\) and \(Z_2\) are \textbf{183 bits each}, and \(k_1\) and \(k_2\) are \textbf{181} and \textbf{183 bits} respectively, while the prime bit lengths are \textbf{1124} and \textbf{1129 bits}. In contrast, RSA uses \textbf{1024-bit} and \textbf{1022-bit} primes without revealing other key parameters. Notably, the execution time of the Enhanced Asymmetric Encryption (\textbf{311.95 ms}) is significantly lower than RSA's execution time (\textbf{1146.63 ms}), demonstrating improved computational efficiency despite the increased key size and cryptographic complexity. 

\begin{figure}[H]
    \centering
    \includegraphics[width=0.775\linewidth]{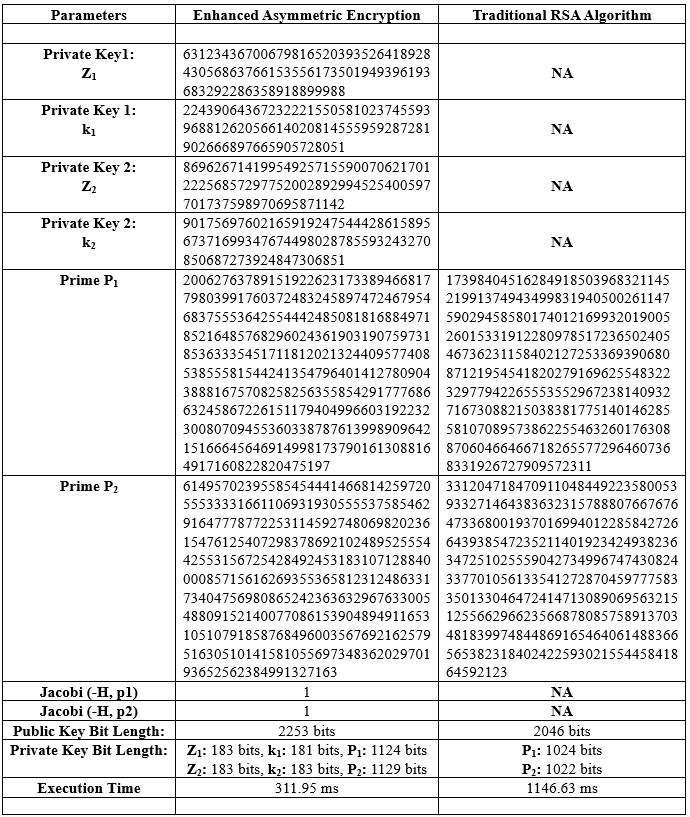}
    \caption{Comparison between Enhanced Asymmetric Encryption and Traditional RSA Method}
    \label{table3}
\end{figure}

\noindent
Besides in cryptography, the polynomial $f_{(Z,k,H)}(n)$ has significant role in other various fields too. Let us discuss another important role of this polynomial in the field of Signal Processing.

\section{Symmetric Prime-Based Channel Allocation for Secure and Fast Signal Transmission}

\noindent Since the polynomial \( f_{(Z,k,H)}(n) \) generates a rich distribution of prime numbers and exhibits symmetrical properties, we can draw a novel analogy by mapping each prime value generated by \( f_{(Z,k,H)}(n) \) to a unique frequency, and each input value \( n \) to a corresponding communication channel. Specifically, due to the symmetric nature of the polynomial for a fixed range \( n \in [0, n_2] \), where \( f_{(Z,k,H)}(n) = f_{(Z,k,H)}(n_2 - n) \), two distinct channel indices (i.e., values of \( n \)) will correspond to the same frequency and \(n_2\) is always an odd positive integer. This implies that two different communication channels can simultaneously operate using the same carrier frequency, yet without interference, due to their mathematical mapping. Traditionally, data are transmitted and received through a single channel at a given frequency \cite{signaltransmission}. However, with this polynomial-based symmetry, both the sender and receiver can exploit the mirrored channels (i.e., \( n \) and \( n_2 - n \)) to transmit data over the same frequency, effectively doubling bandwidth utilization. This not only enhances the speed of data transmission but also introduces a new layer of security, as the pairing of channels through symmetric mapping is non-trivial to predict without knowledge of the polynomial parameters. This framework holds great promise for future communication technologies, especially in the development of highly efficient and secure 6G and 7G systems.

\noindent 
The prime-rich polynomial expressed in eq.(\ref{1.2}) is specifically constructed to produce a high density of prime numbers for integer values of \( n \) over a specified interval. Each output value is assumed to be a prime number and is mapped to a unique carrier frequency. We define:
\[
f_n := f_{(Z,k,H)}(n) \tag{8.1}\label{8.1}
\]
so that the \( n^\text{th} \) communication channel operates at the frequency \( f_n \). The key feature of interest is the symmetry property of this polynomial:
\[
f_{(Z,k,H)}(n) = f_{(Z,k,H)}(n_2 - n) \tag{8.2}\label{8.2}
\]
which holds only for a special choice of the parameter \( Zk \). To ensure symmetry across the range \( [0, n_2] \), we define:
\[
Zk =
\frac{n_2 + 1}{2}, where \quad n_2 \in \{ n \in \mathbb{Z}^+ \mid n \equiv 1 \pmod{2} \} \tag{8.3}\label{8.3}
\]
This choice ensures that the vertex of the quadratic lies at the midpoint of the range and enforces the required symmetry. When eq.(\ref{8.3}) is satisfied, we observe that:
\[
f_n = f_{n_2 - n} \tag{8.4}\label{8.4}
\]
Thus, from eq.(\ref{8.4}) two different channels indexed by \( n \) and \( (n_2 - n) \) are assigned the \emph{same} frequency \( f_n \), but can still transmit distinct data streams. This enables mirrored channel frequency allocation, allowing duplex communication, both transmission and reception, over the same frequency via distinct channels. This approach significantly improves bandwidth utilization and improves security through intrinsic mathematical encoding.

To illustrate this in a real-world setting, consider a futuristic 6G network with two users: Alice and Bob. Let the system allocate channel indices \( n \in [0, 19] \), implying \( n_2 = 19 \). 
\[
\therefore Zk = \frac{19+1}{2} = 10
\]
 Then, the polynomial becomes:
\[
f(n) = n^2 - (2 \cdot 10 - 1)n + \frac{(2 \cdot 10 - 1)^2 + H}{4} = n^2 - 19n + \frac{361 + H}{4}.
\]
As \( H \) belongs to the set of Heegner numbers i.e. 
\(\{1, 2, 3, 7, 11, 19, 43, 67, 163\}\), let us consider \(H=163\),
\[
\therefore f(n) = n^2 - 19n + \frac{524}{4} = n^2 - 19n + 131.
\]
Now suppose that Alice and Bob want to communicate with each other in a more secure and faster way, let us allot a channel \( n = 4 \) to Alice, while Bob uses a channel \( n = 15 \), with \( 4 + 15 = 19 = n_2 \). Then:
\[\boxed{
f(4) = 16 - 76 + 131 = 71, \quad f(15) = 225 - 285+ 131 = 71} \tag{8.5}\label{8.5}
\]
Eq.(\ref{8.5}) shows that on different channels \( n = 4 \) and \( n = 15 \) Alice and Bob can conveniently communicate with each other with the same frequency \(f(4) = f(15) = 71\). This highlights the remarkable symmetry of the polynomial \( f(n) = n^2 - 19n + 131 \), derived using the parameters \( Zk = 10 \) and \( H = 163 \), a Heegner number. In a futuristic 6G communication environment, such symmetry allows the allocation of mirrored channel pairs, such as \( n = 4 \) and \( n = 15 \) that yield identical polynomial values, here \( f(4) = f(15) = 71 \). This intrinsic symmetry can be used to optimize channel pairing, error detection, or even secure transmission strategies. As a result, the communication infrastructure of the future (notably 6G and 7G) can leverage this mathematical structure to double channel availability, reduce frequency contention, and embed natural cryptographic security through prime-based frequency encoding.

\section{Conclusion and Future Scope}

In this paper, we have investigated the prime-rich nature of the polynomial family \( f_{(Z,k,H)}(n) = n^2 - (2Zk - 1)n + \frac{(2Zk - 1)^2 + H}{4} \). Through rigorous mathematical analysis and numerical experiments, we have demonstrated the ability of this family to consistently generate a higher density of primes within a defined range of \( n \). By employing deep theoretical frameworks such as the Prime Number Theorem, the Bateman-Horn Conjecture, the Riemann Hypothesis, and the properties of the Riemann Zeta function, we validated the prime-generating efficiency and theoretical consistency of our approach. Moreover, we presented an approximation heuristic for the Bateman-Horn constant \( C_f \) based on the quadratic residue asymmetry, which aligns remarkably well with the actual values for carefully chosen parameters, thus reinforcing the robustness of our method.
Prime-rich polynomials such as \( f_{(Z,k,H)}(n) \) can play a vital role in computational number theory, cryptographic key generation, and secure digital communications. Given their prime-density property, these polynomials offer potential improvements in pseudorandom prime generation algorithms, primality testing frameworks, and lightweight cryptographic systems for constrained environments.
Further research could involve extending this polynomial family to higher degrees or multivariate settings, studying their behavior under modular arithmetic and residue classes, and analyzing their distribution patterns in the context of algebraic number theory. Exploring connections with elliptic curves, L-functions, and modular forms could also yield deeper insights. Furthermore, investigating machine learning models trained on prime-rich polynomial data could provide new ways to predict and classify prime distributions more efficiently. 

\bibliographystyle{alpha}
\bibliography{mybib}

\end{document}